\documentclass[aps,showkeys,showpacs,notitlepage,pra,twocolumn]{revtex4-1}
\usepackage{graphicx} 
 \usepackage{amsmath}

\usepackage[pdfauthor={Richard J. Mathar}, colorlinks=true,pdfkeywords={2-regular digraphs}]{hyperref}

\begin{document}

\title{2-regular Digraphs of the Lovelock Lagrangian}

\author{Richard J. Mathar}
\homepage{http://www.mpia-hd.mpg.de/~mathar} 
\affiliation{Max-Planck Institute of Astronomy, K\"onigstuhl 17, 69117 Heidelberg, Germany} 

\date{\today}

\begin{abstract}
The manuscripts tabulates arc lists of the 1, 1, 3, 8, 25, 85, 397 \ldots unlabeled 2-regular digraphs on 
$n=0, 1, 2, \ldots,9$
nodes, including disconnected graphs, graphs with multiarcs and/or graphs with loops.
Each of these graphs represents one term of the Lagrangian of Lovelock's type ---
a contraction of a product of $n$ Riemann tensors --- once the 2 covariant and 2 contravariant
indices of a tensor are associated with the in-edges and out-edges of a node.
\end{abstract}

\maketitle

\section{Regular Digraphs}
\subsection{Nomenclature}
In the common language of graph theory, digraphs (directed graphs) are graphs where
the edges (called arcs) are oriented, i.e., the two nodes of an arc are distinguished to be tail and head
of the arc. The indegree of a node is the number of arcs that
have heads at a node; the outdegree is the number of arcs that have tails at a node.
$k$-regular digraphs are digraphs where the indegree and outdegree at each node is the
same $k$ \cite{BrinkmannDM313,RamanathJGT11}.
Loops are arcs where head and tail are the same node. (Each loop increases the
indegree and outdegree of the node by one.)
Multiarcs are multisets of two or more arcs that share a head and tail node.

The underlying simple graph of a digraph is obtained by reducing the arcs to undirected edges,
replacing multiedges by single edges, and removing loops. We will be only interested
in these simple graphs to classify the digraphs by the number of components of their
underlying simple graphs. (That means according to \textit{weak} connectivity).

(Vertex) labeled graphs have distinct labels at their nodes, usually taken
to be the positive integers while handling graphs on computers, or letters from \textit{a} onwards.
We define the Adjacency Matrix of a digraph of nodes labeled $1,2,\ldots$ as the square
matrix which contains in row $r$ and column $c$ the number of arcs which start at node $r$ 
and end at node $c$. Construction of the labeled $k$-regular digraphs admitting multiarcs and loops
is therefore equivalent
to constructing the $n\times n$ matrices with non-negative integer entries where all
row sums and all column sums are $k$.
These are registered in Table A008300 in the Online Encyclopedia of Integer Sequences
for the cases where multiarcs are \emph{not} admitted, and in Table A257493
if multiarcs are admitted \cite{sloane}.

Denote the number of unlabeled and labeled $k$-regular digraphs with $n$ nodes
and $c$ components by $U_k(n,c)$ and $L_k(n,c)$. 

The unlabeled graphs with $c>1$ components can be derived from
the connected graphs quickly
by the Multiset Transformation, collecting arrangements
over all partitions of $n$ into $c$ parts \cite[Theor. I.1]{Flajolet}\cite[(4.2.3)]{Harary}:
\begin{equation}
U_k(n,c) = \sum_{\stackrel{n=n_1+2n_2+\cdots +cn_c}{n_i\ge 0}} 
\prod_{i=1}^c \binom{U_k(i,1)+n_i-1}{n_i}.
\end{equation}

The number of unlabeled and labeled
$k$-regular digraphs with $n$ nodes 
is
\begin{eqnarray}
U_k(0)=1;\quad U_k(n) & \equiv & \sum_{c=1}^n U_k(n,c); \\
L_k(0)=1;\quad L_k(n) & \equiv & \sum_{c=1}^n L_k(n,c).
\end{eqnarray}
The case of the unlabeled 1-regular digraphs is simple: the connected
1-regular digraphs are cycles (unicycles), including the case with 1 node and its loop: 
\begin{equation}
U_1(n,1) =1.
\end{equation}
So the
number of unlabeled, not necessarily connected, 1-regular digraphs on $n$ nodes, $U_1(n)$,
is the number of partitions of $n$ \cite[A000041]{sloane},
and the $U_1(n,c)$ are the partition numbers \cite[A008284]{sloane}.

\subsection{Symmetries}

The key part of this work is to identify the Automorphism Group of the labeled
2-regular digraphs, which is the group of permutations of the labels which keeps
the Adjacency Matrix of a graph the same. This bundles a set of one or more
labeled digraphs which are mapped onto each other by the permutations of
the group, and each such set is represented by a single unlabeled 2-regular digraph.
This is one way of erasing/forgetting any particular order
on the nodes while maintaining the connectivity and structure of the graph.

The unlabeled 2-regular digraphs might also be obtained by starting from the cubic QED 
vacuum polarization diagrams 
\cite{MelloPRD85,MatharVixra1901},
coalescing each edge that represents an undirected photon line and its two incident nodes
into a single node such that only the directed fermion lines connect to the remaining nodes---representing
nodes of a $\varphi^4$ theory---, 
plus a cleanup that eliminates duplicates.
Another approach is to construct first the 4-regular undirected unlabeled graphs 
\cite{MelloPRD85}\cite[A129429]{sloane}, and to add directions later.

We shall keep track of the Automorphism Group $\cal A$ of each graph by writing down the
Cycle Index of the group of the node permutations \cite[\S 2]{Harary}. Since $\cal A$
is  a subgroup of the group of all permutations of $n$, its order, $|\cal A|$, is a divisor of $n!$.
The denominator of the cycle index polynomial is $|\cal A|$, so we can recover
the number of labeled digraphs by dividing $n!$ through the denominator of the Cycle Index \cite[(1.1.3)]{Harary}.
For $n=3$ nodes, for example, $U_2(3)=8$ unlabeled 2-regular digraphs exist
with 4 different Cycle Indices:
\begin{itemize}
\item
1 graph with cycle index $(t_1^3)/1$,
\item
3 graphs with cycle index $(t_1^3+t_1t_2)/2$,
\item
2 graphs with cycle index $(t_1^3+2t_3)/3$,
\item
and 2 graphs with cycle index $(t_1^3+3t_1t_2+2t_3)/6$,
\end{itemize}
and $1\times 3!/1 + 3\times 3!/2 + 2\times 3!/3 + 2\times 3!/6=21 = L_2(3)$ 
is the number of labeled
graphs on 3 nodes \cite[A000681]{sloane}.
As we are admitting loops, there is always at least one graph on $n$ nodes (the one consisting of 
$n$ isolated nodes with two loops each) that has the maximum symmetry here, $|{\cal A}|=n!$.

The labeled graphs with $c>1$ components
are deduced from the weakly connected labeled graphs by a Bell transformation,
summing over all compositions of $n$ into positive parts weighted by multinomial
coefficients \cite[EIJ]{BowerOeisT2}:
\begin{equation}
L_k(n,c) = \frac{1}{c!}\sum_{\stackrel{n=n_1+n_2+\cdots +n_c}{n_i\ge 1}} 
\binom{n}{n_1,n_2,\cdots, n_c}
\prod_{i=1}^c L_k(n_i,1).
\end{equation}
If
\begin{equation}
L_k(x,1) \equiv \sum_{n\ge 1} \frac{L_k(n,1)}{n!}x^n
\end{equation}
denotes the exponential generating function of the weakly connected labeled graphs,
the bivariate exponential generating function of the labeled graphs is \cite{GilbertCDM8}
\begin{equation}
L_k(x,t) \equiv \sum_{n,c\ge 0}\frac{L_k(n,c)}{n!}x^n t^c
= \exp\left[ t L_k(x,1)\right].
\end{equation}

The corresponding derivation starting from $L_1(n,1)=1$ demonstrates that $L_1(n,c)$
are the Stirling Numbers of the Second Kind \cite[A008277]{sloane}
and $L_1(n)$ the Bell Numbers \cite[A000110]{sloane}.

For $L_2$, the number of labeled digraphs refined by the number of weakly connected components
is summarized in Table \ref{tab.L}.  

\begin{table}[htb]
\begin{ruledtabular}
\begin{tabular}{r|rrrrrrr|r}
$n\backslash c$ & $1$ & $2$& 3 & 4 & 5 & 6 & 7 &  $L_2(n)$ \\
\hline
1 & 1 &&&&&&& 1\\
2 & 2 & 1 &&&&&& 3\\
3 & 14 & 6 & 1 &&&&& 21 \\
4 & 201 & 68 & 12 & 1 &&&& 282 \\
5 & 4704 & 1285 & 200 & 20 & 1 &&& 6210\\
6 & 160890 & 36214 & 4815 & 460 & 30 & 1 && 202410\\
7 & 7538040 & 1422288 & 160594 & 13755 & 910 & 42 & 1 & 9135630 \\
\end{tabular}
\end{ruledtabular}
\caption{Labeled 2-regular digraphs $L_2(n,c)$ with $n$ nodes and $c$ (weak) components, 
allowing loops and multiarcs 
\cite[A307804]{sloane}.
} 
\label{tab.L}
\end{table}

\begin{widetext}

\section{Gallery of Unlabeled 2-regular Digraphs}
The unlabeled 2-regular digraphs on $n\le 5$ nodes are illustrated in the following
sections. Each digraph has $2n$ arcs.
For easier visual recognition, the graphs with more than one component are surrounded
by a frame. 
These visualize throw schemes  for
$n$ 2-armed jugglers.
\subsection{1 graph on 1 node}\label{sec.1}

\includegraphics[scale=0.3]{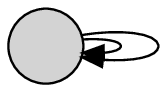}

\subsection{3 graphs (2 connected) on 2 nodes}

\includegraphics[scale=0.3]{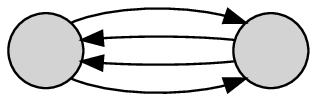}
\includegraphics[scale=0.3]{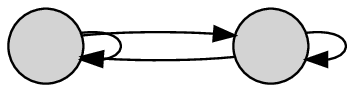}
\fbox{\includegraphics[scale=0.3]{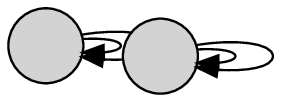}}

\subsection{8 graphs (5 connected) on 3 nodes}

\includegraphics[scale=0.3]{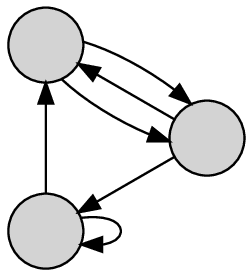}
\includegraphics[scale=0.3]{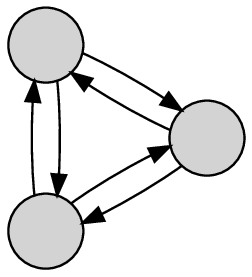}
\includegraphics[scale=0.3]{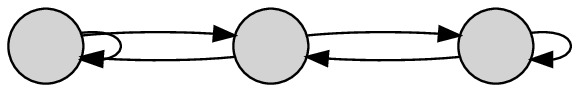}
\includegraphics[scale=0.3]{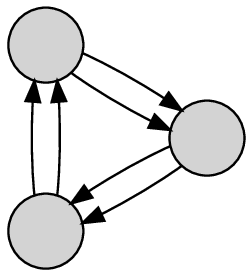}

\fbox{\includegraphics[scale=0.3]{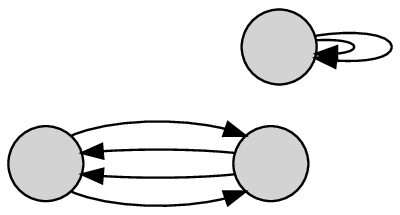}}
\includegraphics[scale=0.3]{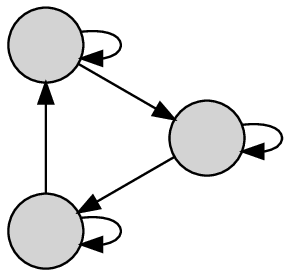}
\fbox{\includegraphics[scale=0.3]{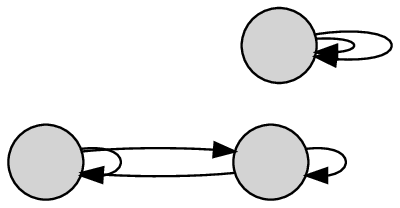}}
\fbox{\includegraphics[scale=0.3]{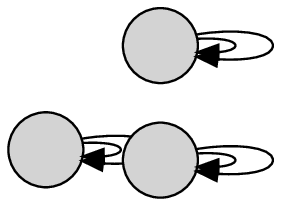}}

\subsection{25 graphs (14 connected) on 4 nodes}

\includegraphics[scale=0.3]{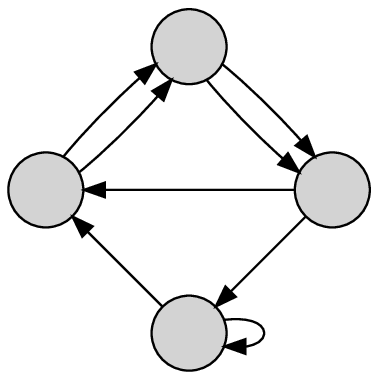}
\includegraphics[scale=0.3]{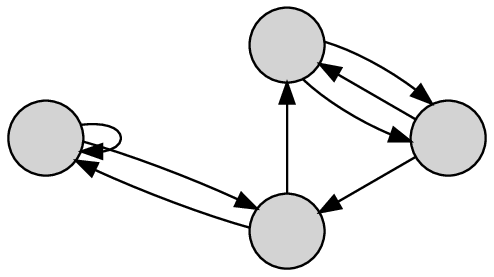}
\includegraphics[scale=0.3]{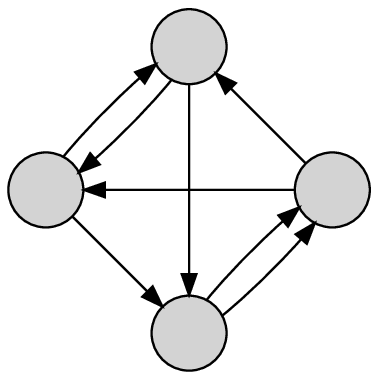}
\includegraphics[scale=0.3]{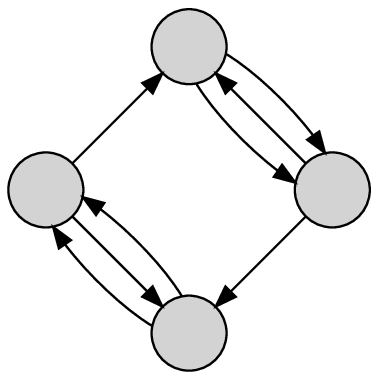}

\includegraphics[scale=0.3]{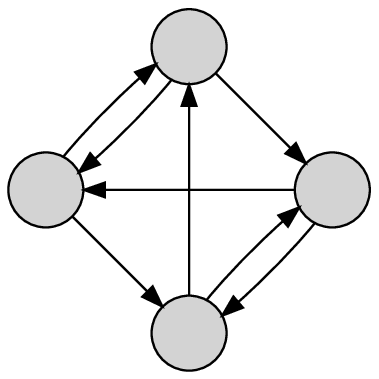}
\includegraphics[scale=0.3]{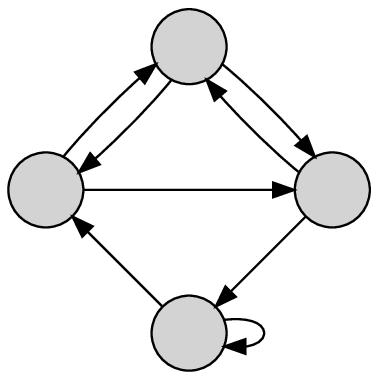}
\includegraphics[scale=0.3]{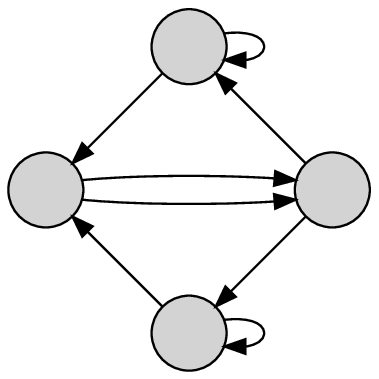}
\includegraphics[scale=0.3]{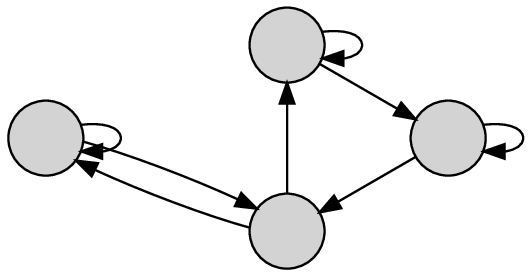}

\includegraphics[scale=0.3]{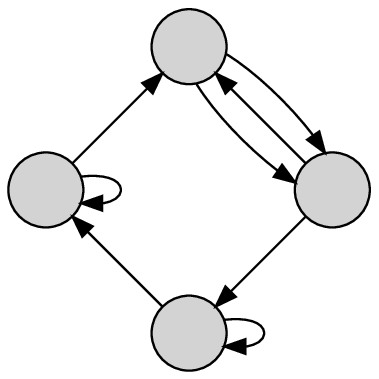}
\includegraphics[scale=0.3]{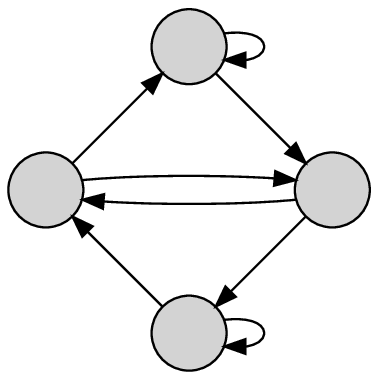}
\fbox{\includegraphics[scale=0.3]{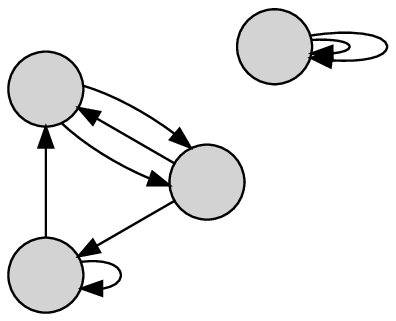}}
\includegraphics[scale=0.3]{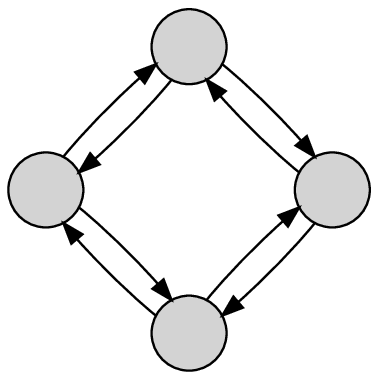}

\includegraphics[scale=0.3]{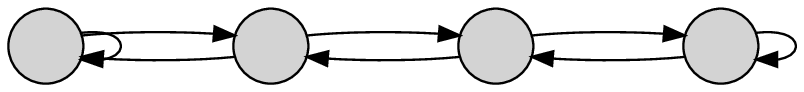}
\fbox{\includegraphics[scale=0.3]{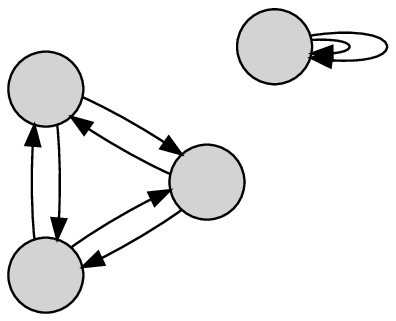}}
\fbox{\includegraphics[scale=0.3]{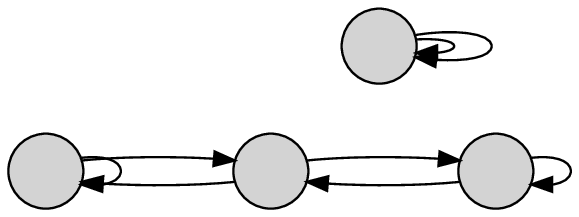}}
\fbox{\includegraphics[scale=0.3]{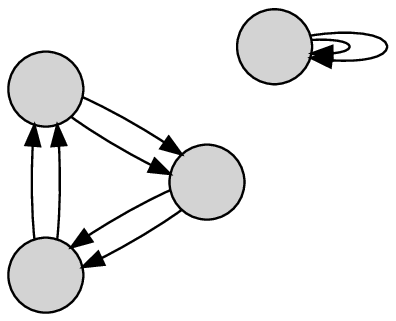}}

\includegraphics[scale=0.3]{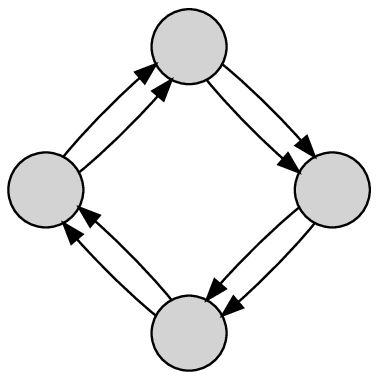}
\fbox{\includegraphics[scale=0.3]{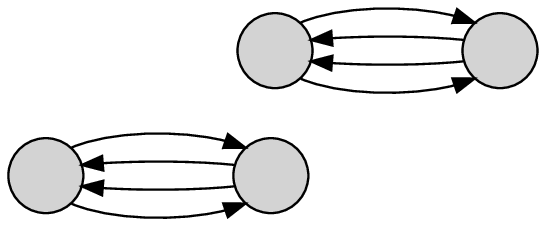}}
\fbox{\includegraphics[scale=0.3]{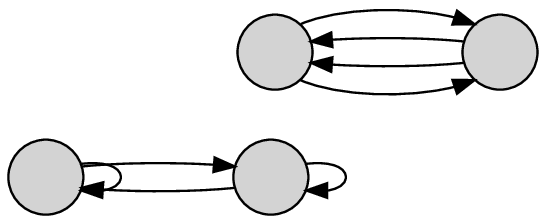}}

\fbox{\includegraphics[scale=0.3]{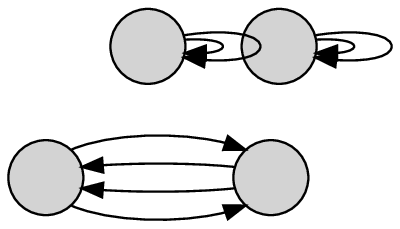}}
\fbox{\includegraphics[scale=0.3]{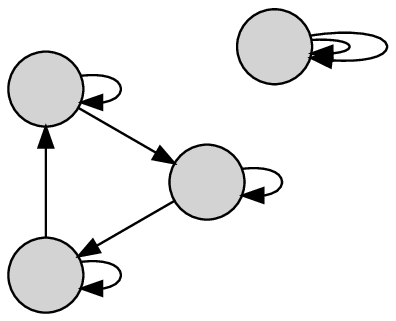}}
\includegraphics[scale=0.3]{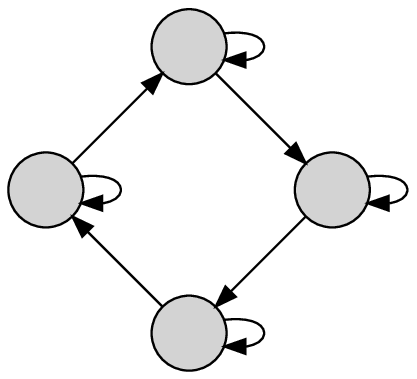}

\fbox{\includegraphics[scale=0.3]{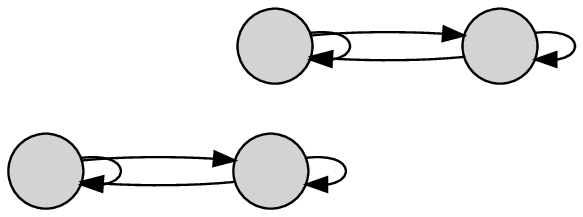}}
\fbox{\includegraphics[scale=0.3]{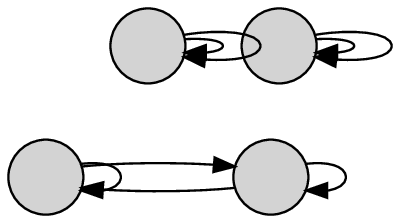}}
\fbox{\includegraphics[scale=0.3]{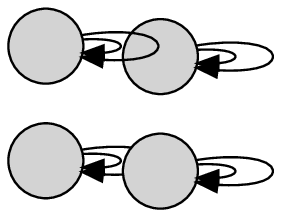}}

\subsection{85 graphs (50 connected) on 5 nodes}

\includegraphics[scale=0.3]{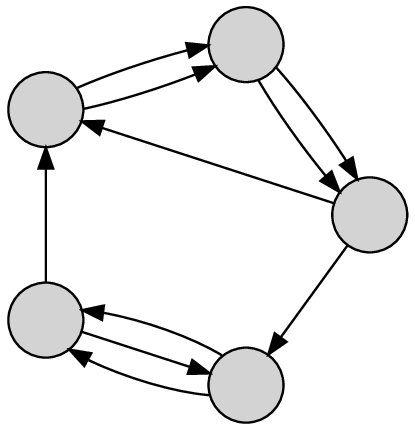}
\includegraphics[scale=0.3]{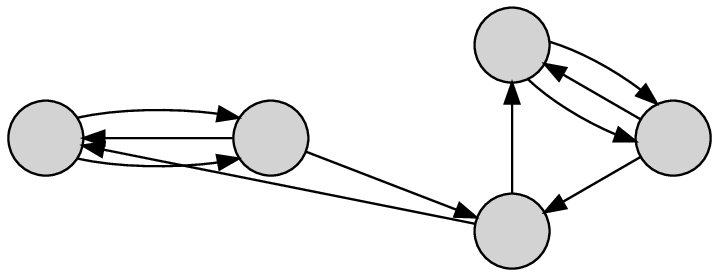}
\includegraphics[scale=0.3]{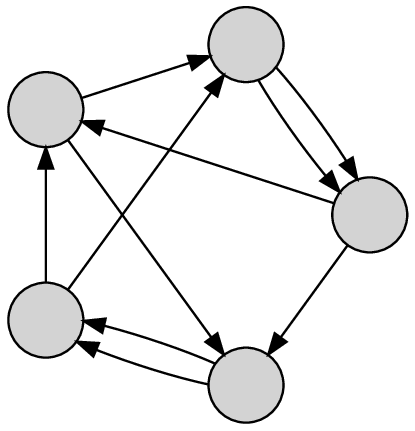}
\includegraphics[scale=0.3]{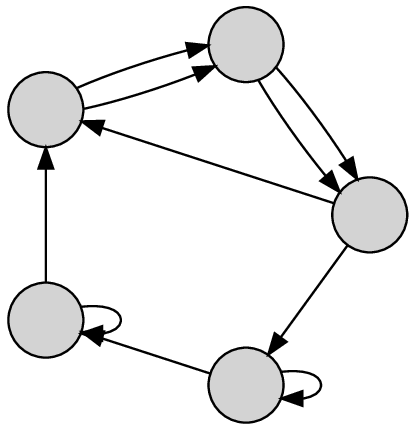}

\includegraphics[scale=0.3]{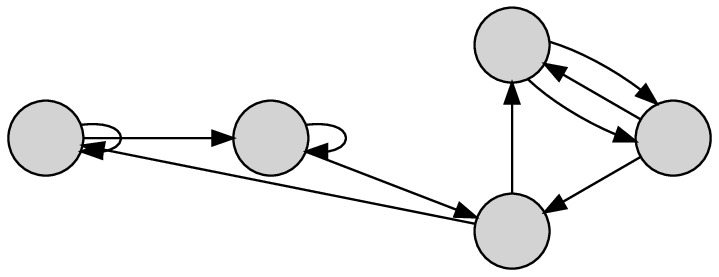}
\includegraphics[scale=0.3]{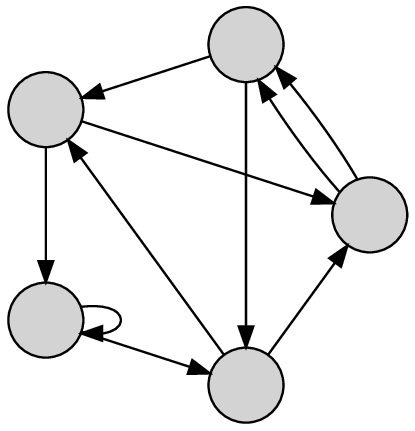}
\includegraphics[scale=0.3]{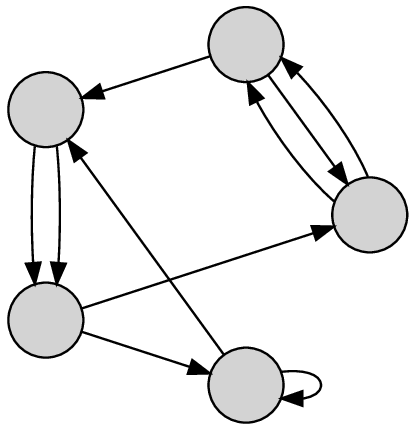}
\includegraphics[scale=0.3]{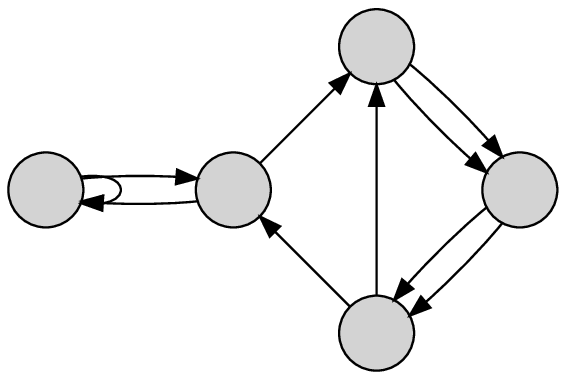}

\includegraphics[scale=0.3]{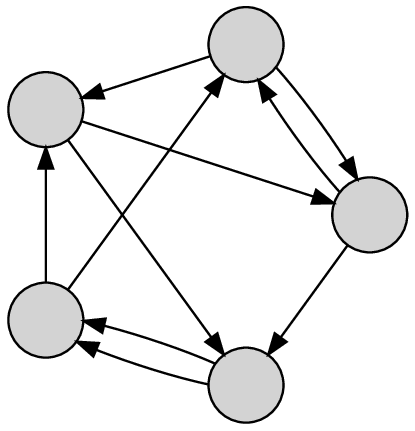}
\fbox{\includegraphics[scale=0.3]{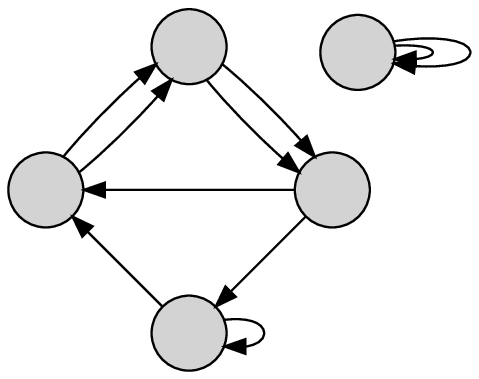}}
\includegraphics[scale=0.3]{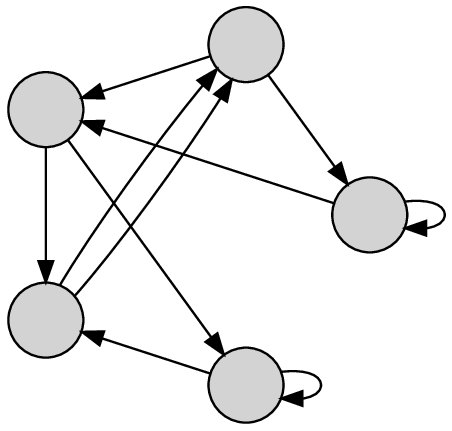}
\includegraphics[scale=0.3]{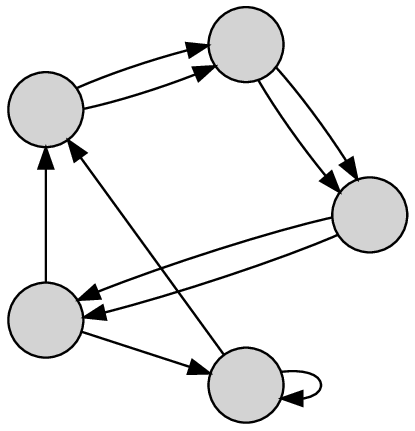}

\includegraphics[scale=0.3]{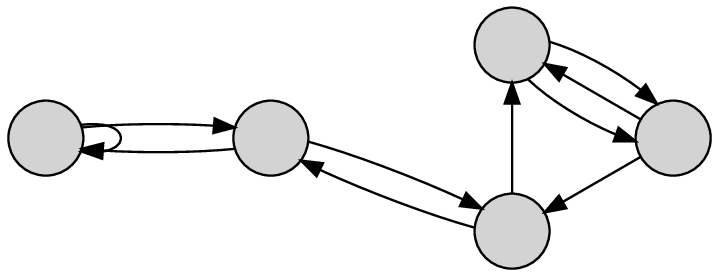}
\includegraphics[scale=0.3]{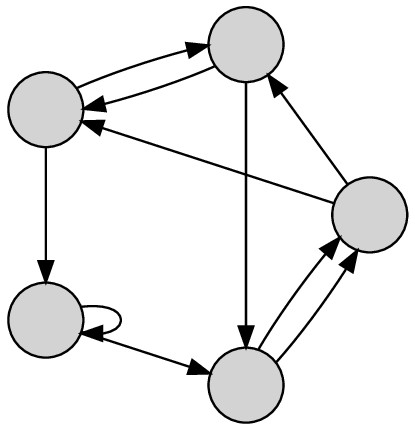}
\includegraphics[scale=0.3]{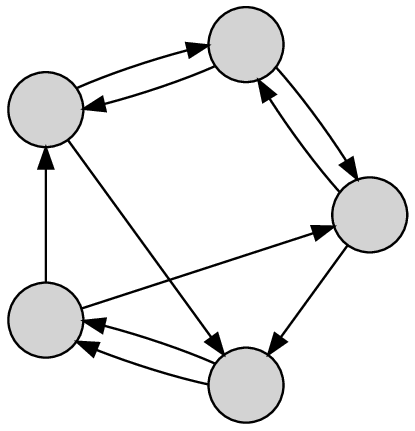}

\includegraphics[scale=0.3]{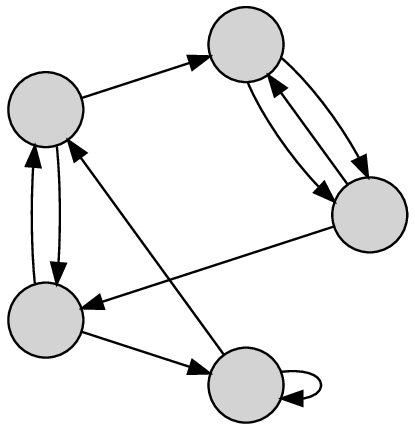}
\includegraphics[scale=0.3]{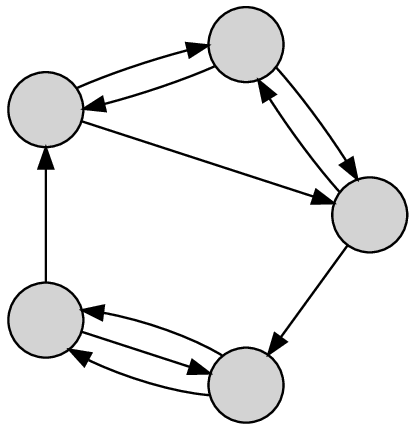}
\fbox{\includegraphics[scale=0.3]{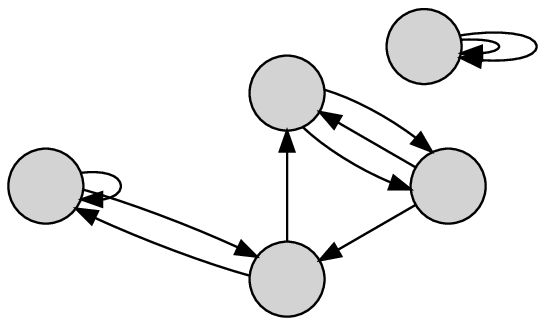}}
\includegraphics[scale=0.3]{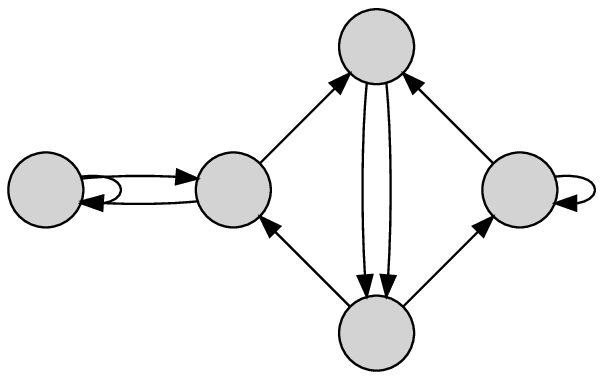}

\includegraphics[scale=0.3]{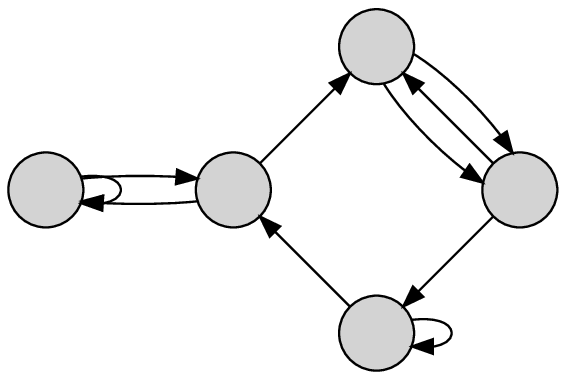}
\fbox{\includegraphics[scale=0.3]{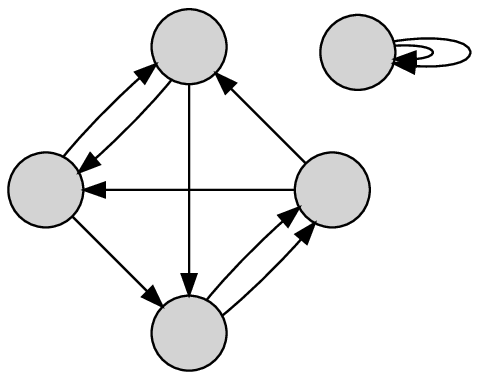}}
\includegraphics[scale=0.3]{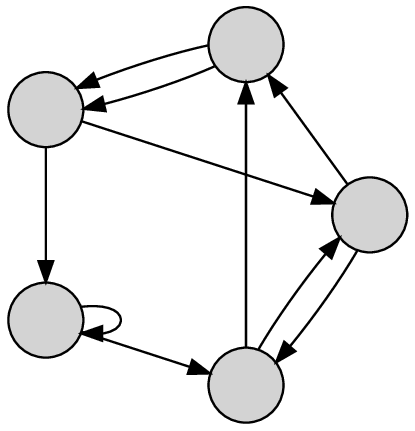}
\includegraphics[scale=0.3]{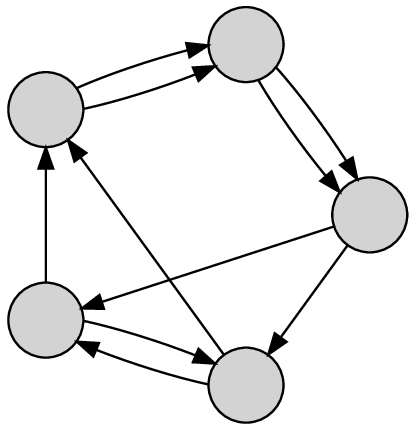}

\fbox{\includegraphics[scale=0.3]{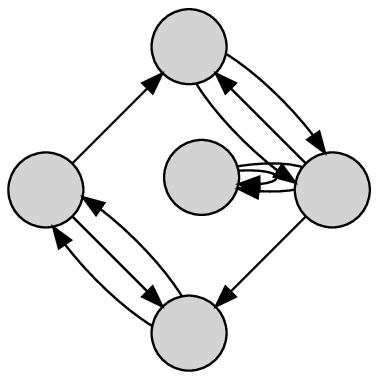}}
\includegraphics[scale=0.3]{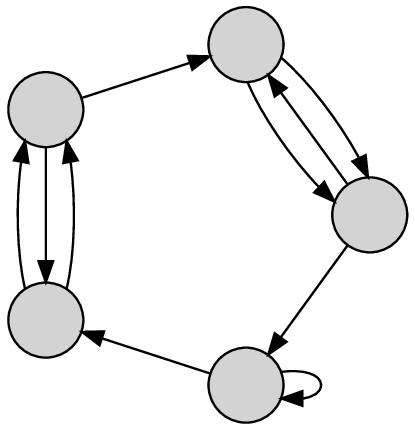}
\includegraphics[scale=0.3]{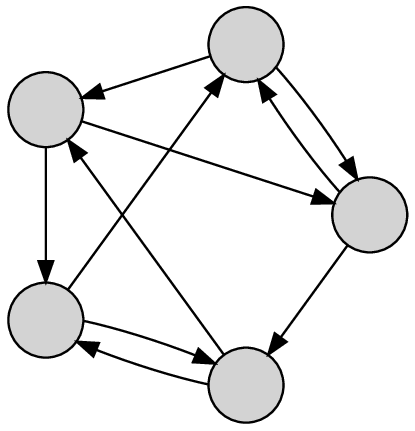}
\includegraphics[scale=0.3]{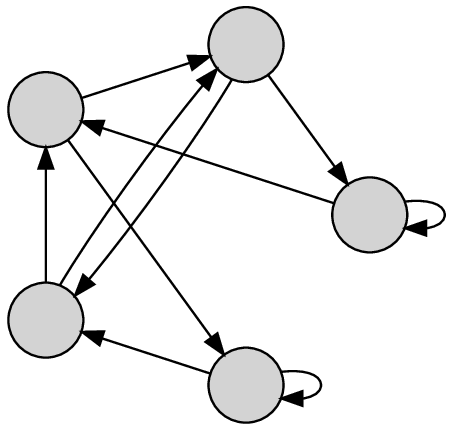}

\includegraphics[scale=0.3]{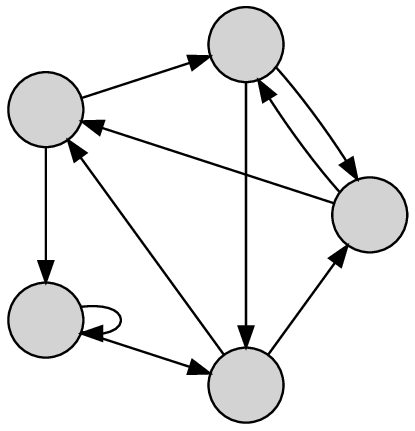}
\includegraphics[scale=0.3]{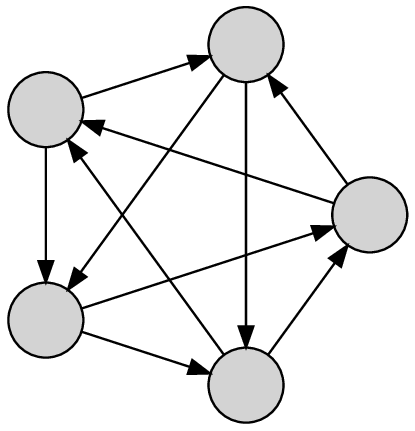}
\includegraphics[scale=0.3]{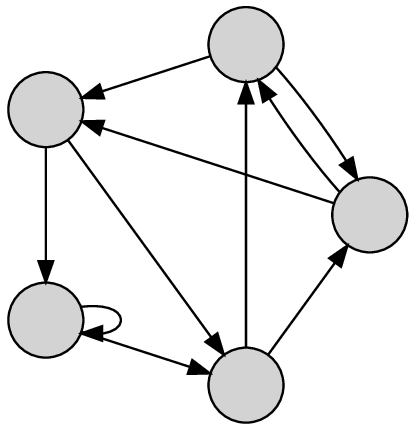}
\includegraphics[scale=0.3]{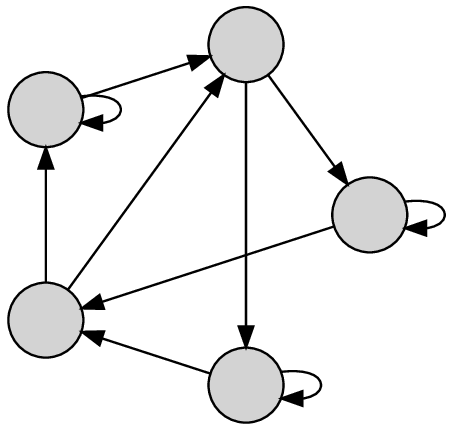}

\includegraphics[scale=0.3]{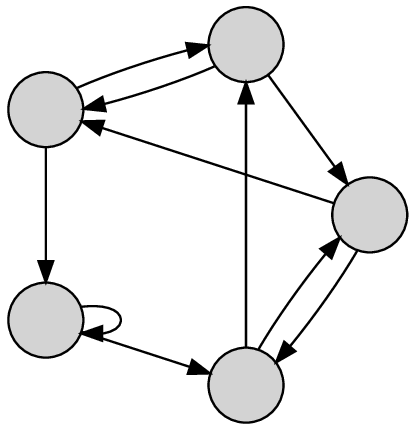}
\includegraphics[scale=0.3]{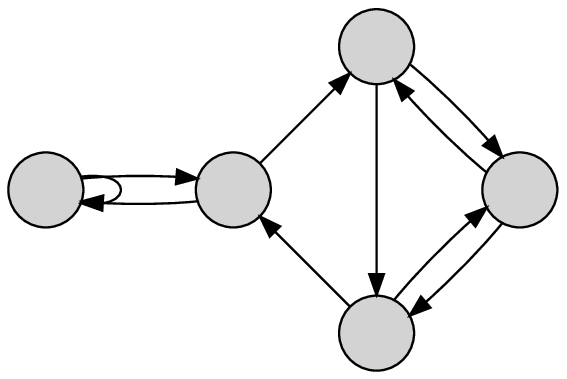}
\includegraphics[scale=0.3]{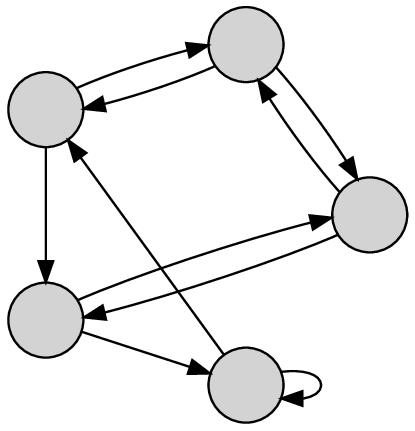}
\includegraphics[scale=0.3]{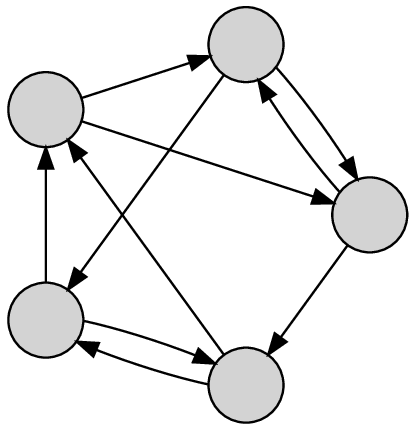}

\includegraphics[scale=0.3]{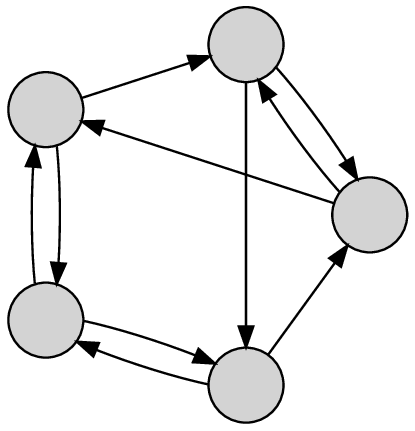}
\includegraphics[scale=0.3]{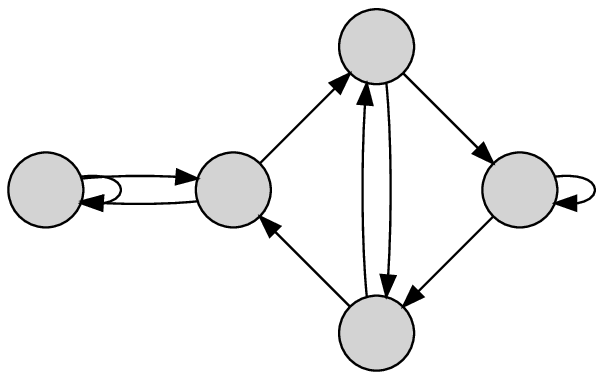}
\includegraphics[scale=0.3]{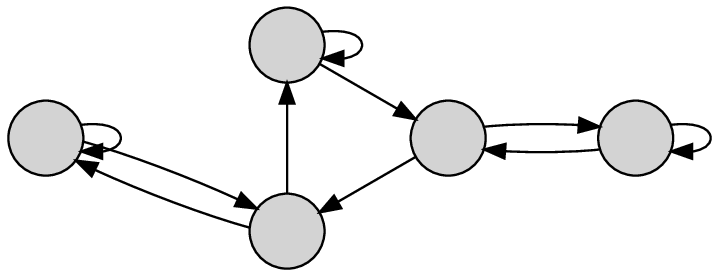}
\includegraphics[scale=0.3]{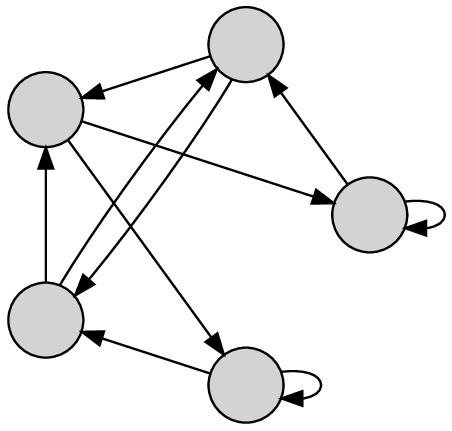}

\includegraphics[scale=0.3]{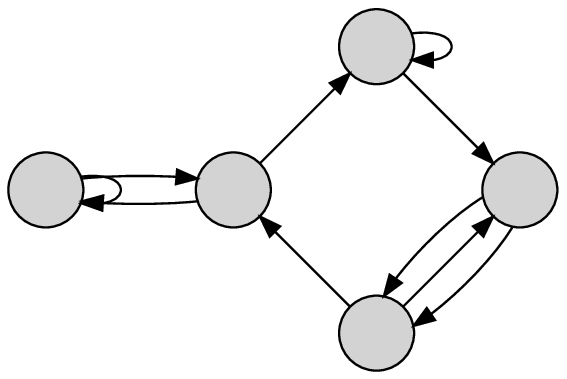}
\includegraphics[scale=0.3]{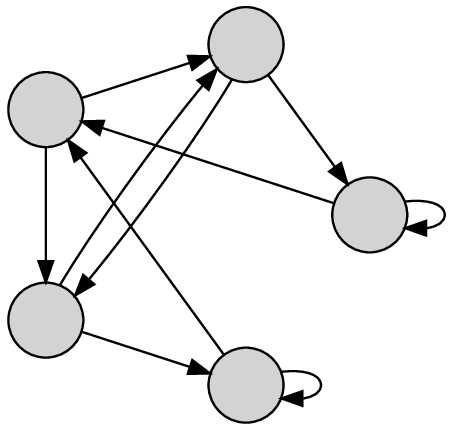}
\fbox{\includegraphics[scale=0.3]{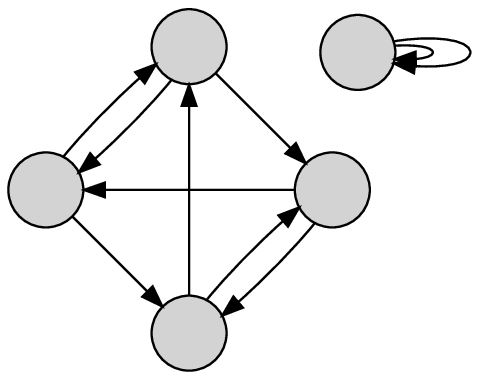}}
\fbox{\includegraphics[scale=0.3]{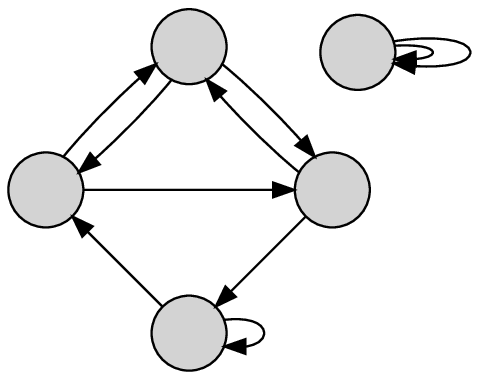}}

\includegraphics[scale=0.3]{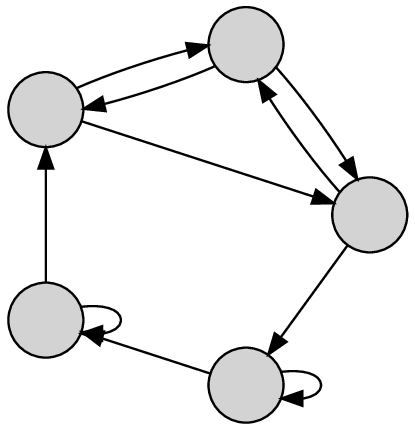}
\fbox{\includegraphics[scale=0.3]{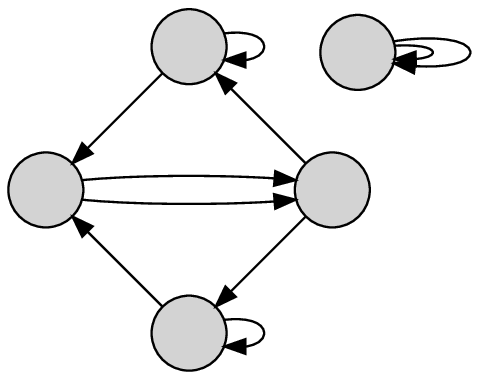}}
\includegraphics[scale=0.3]{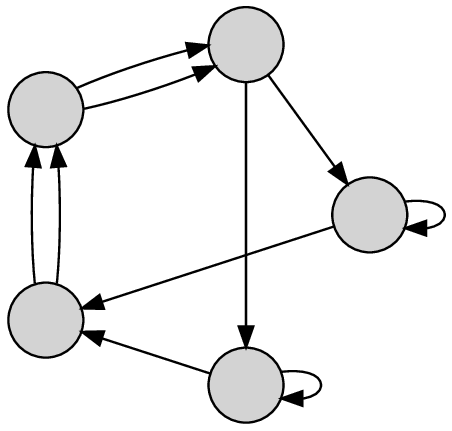}
\includegraphics[scale=0.3]{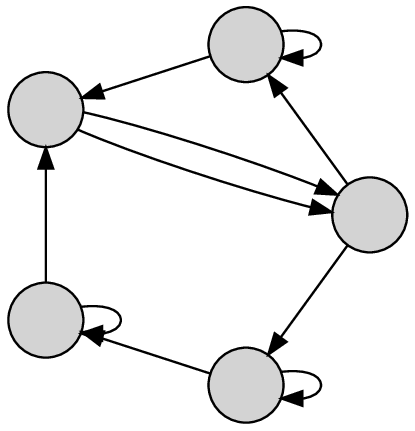}

\includegraphics[scale=0.3]{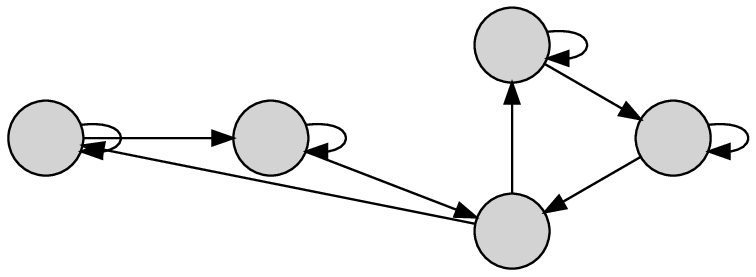}
\includegraphics[scale=0.3]{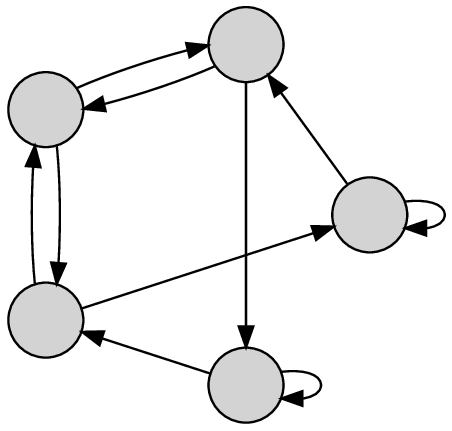}
\includegraphics[scale=0.3]{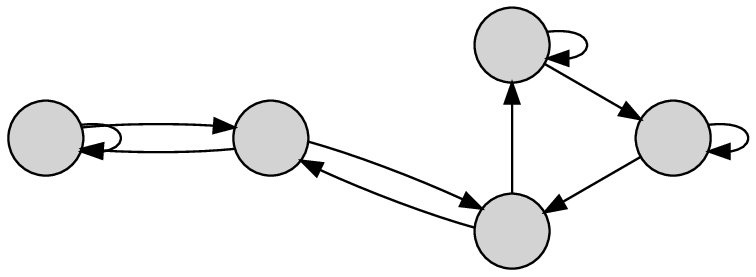}
\includegraphics[scale=0.3]{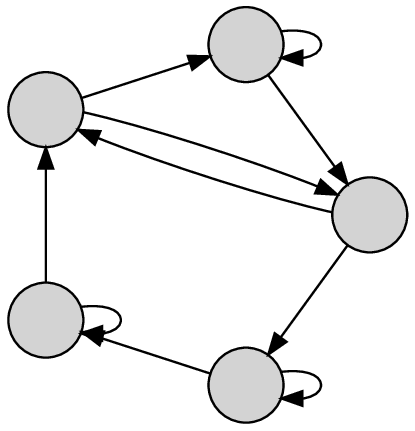}

\fbox{\includegraphics[scale=0.3]{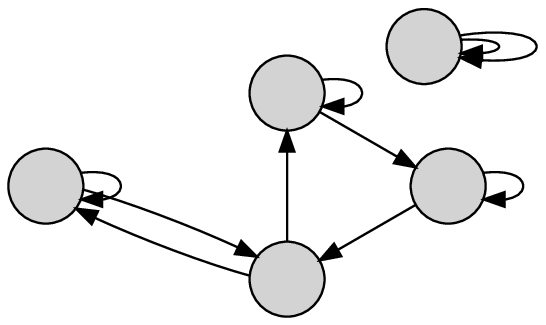}}
\includegraphics[scale=0.3]{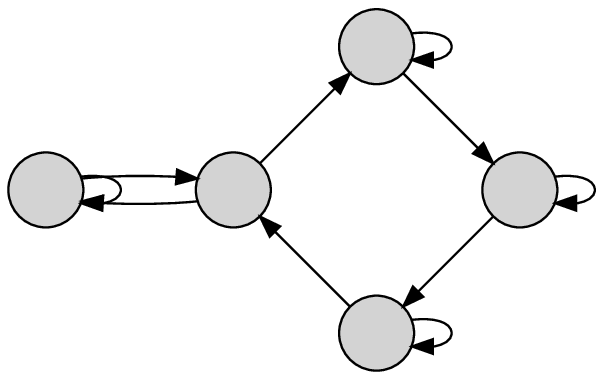}
\fbox{\includegraphics[scale=0.3]{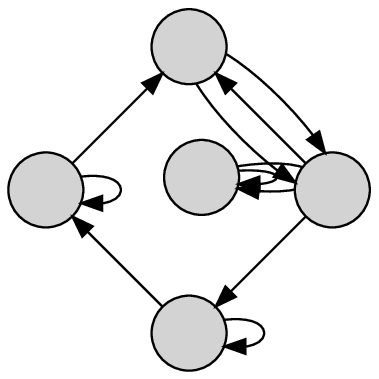}}
\includegraphics[scale=0.3]{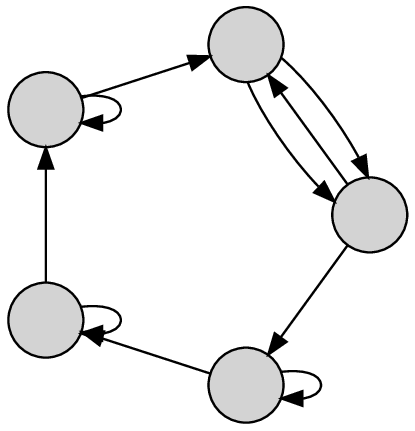}

\fbox{\includegraphics[scale=0.3]{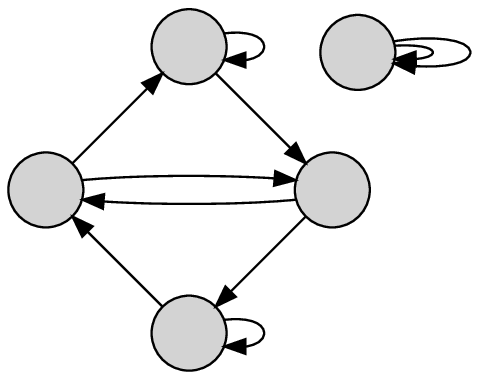}}
\fbox{\includegraphics[scale=0.3]{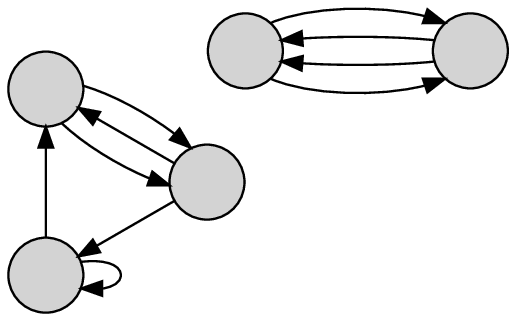}}
\fbox{\includegraphics[scale=0.3]{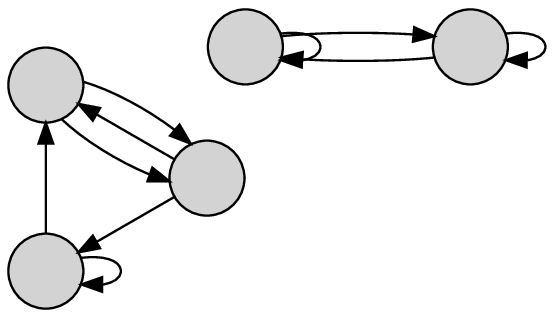}}
\fbox{\includegraphics[scale=0.3]{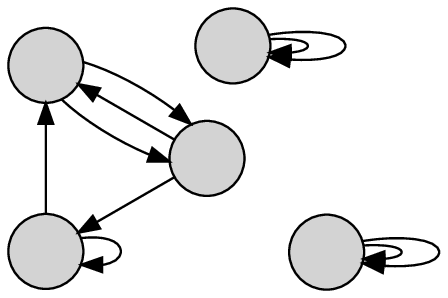}}

\fbox{\includegraphics[scale=0.3]{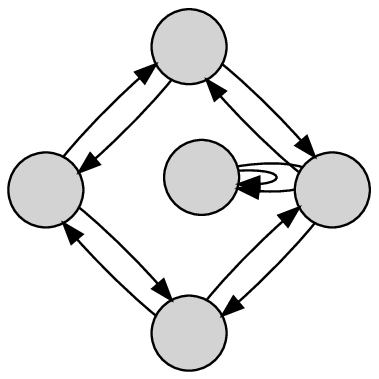}}
\includegraphics[scale=0.3]{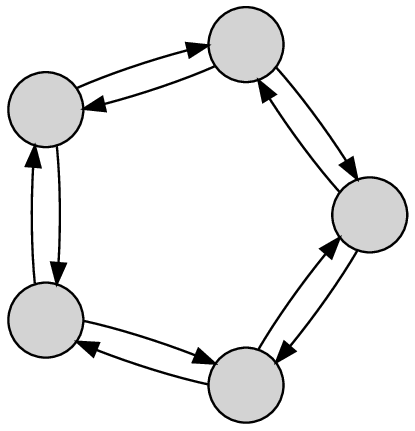}
\includegraphics[scale=0.3]{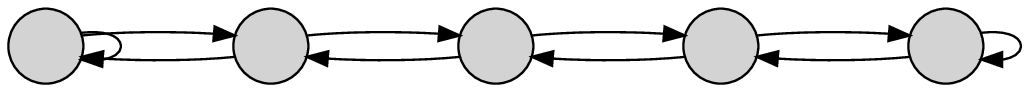}

\fbox{\includegraphics[scale=0.3]{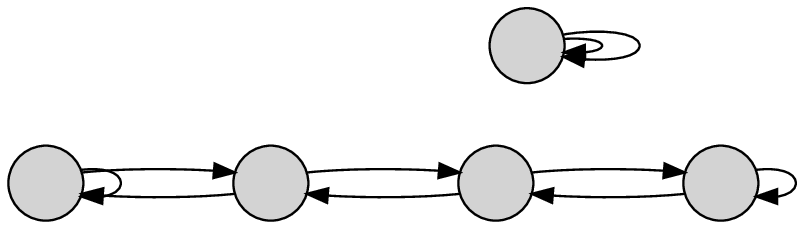}}
\fbox{\includegraphics[scale=0.3]{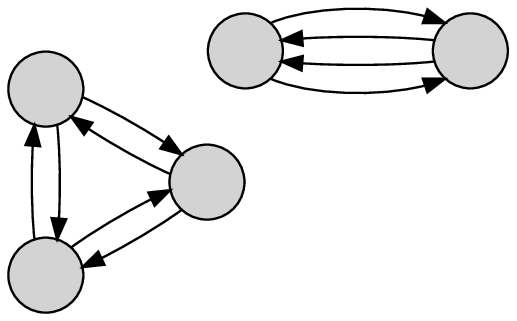}}
\fbox{\includegraphics[scale=0.3]{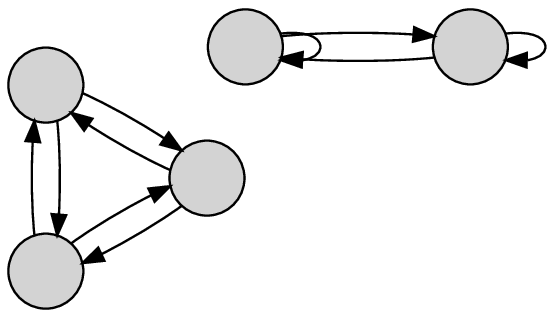}}
\fbox{\includegraphics[scale=0.3]{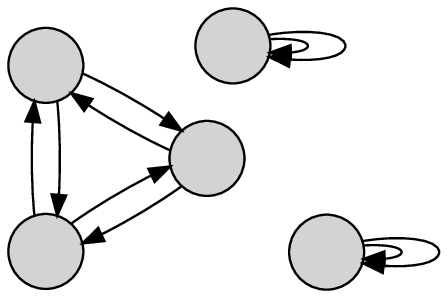}}

\fbox{\includegraphics[scale=0.3]{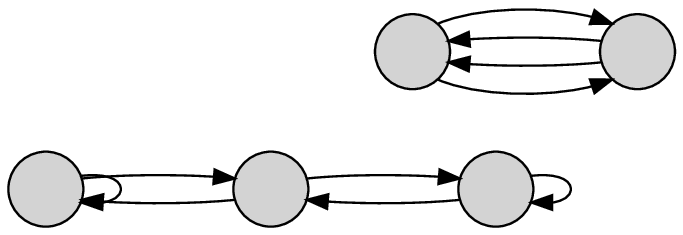}}
\fbox{\includegraphics[scale=0.3]{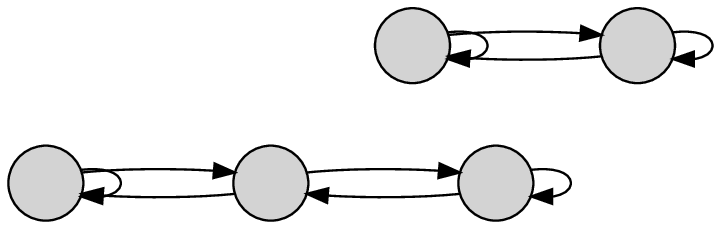}}
\fbox{\includegraphics[scale=0.3]{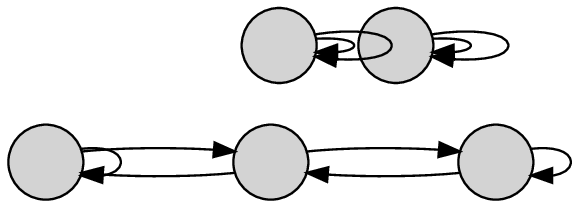}}

\fbox{\includegraphics[scale=0.3]{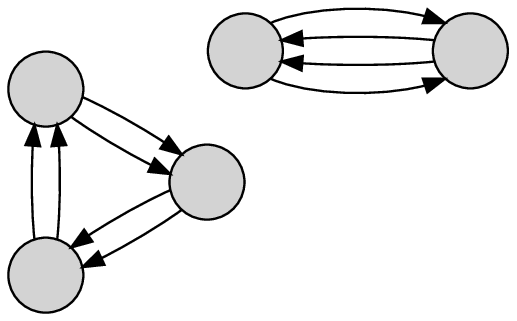}}
\includegraphics[scale=0.3]{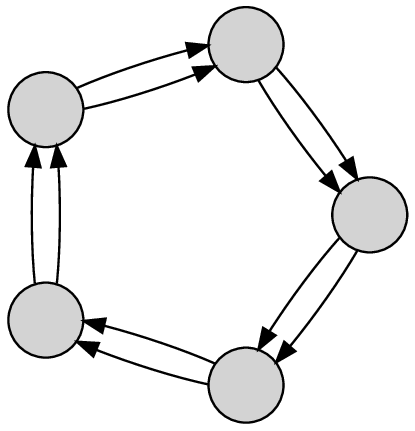}
\fbox{\includegraphics[scale=0.3]{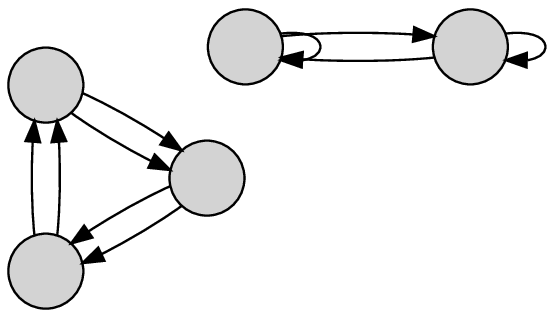}}
\fbox{\includegraphics[scale=0.3]{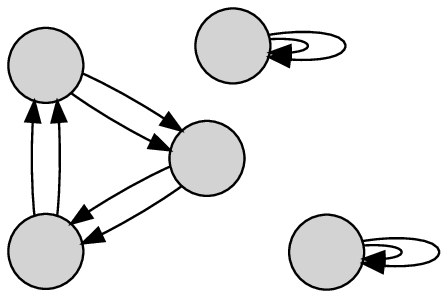}}

\fbox{\includegraphics[scale=0.3]{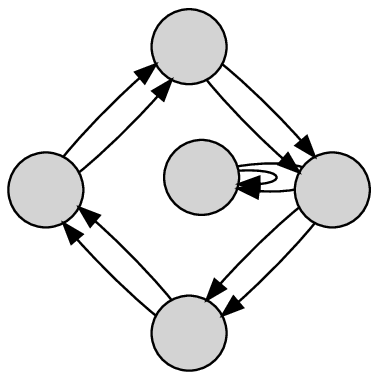}}
\fbox{\includegraphics[scale=0.3]{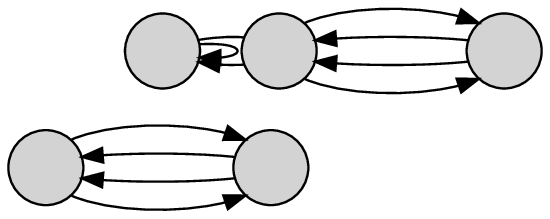}}
\fbox{\includegraphics[scale=0.3]{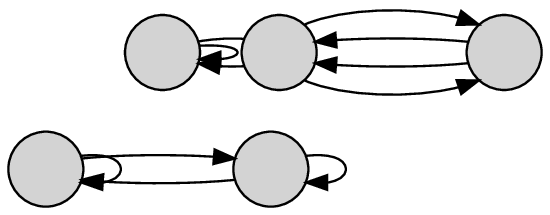}}
\fbox{\includegraphics[scale=0.3]{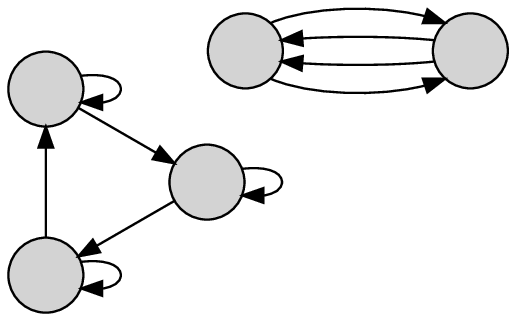}}

\fbox{\includegraphics[scale=0.3]{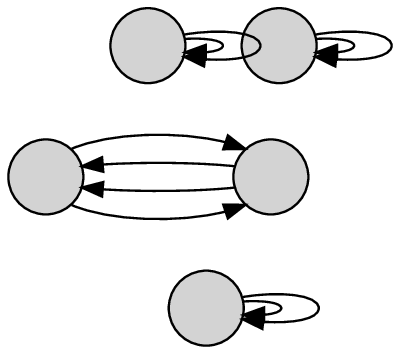}}
\fbox{\includegraphics[scale=0.3]{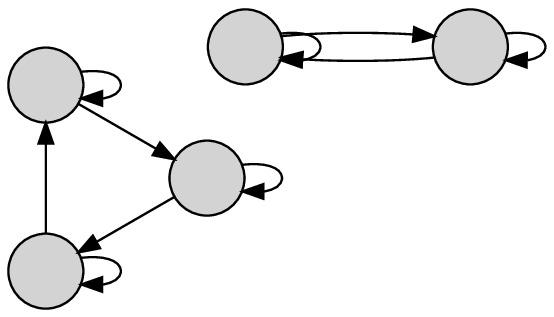}}
\includegraphics[scale=0.3]{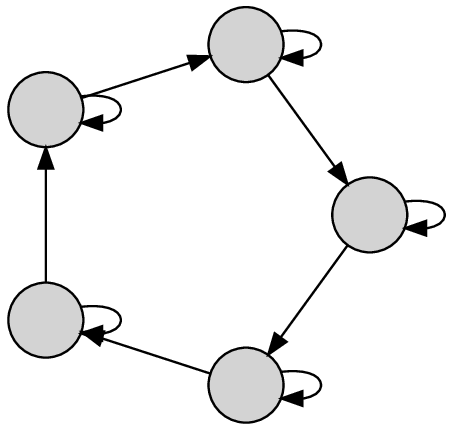}
\fbox{\includegraphics[scale=0.3]{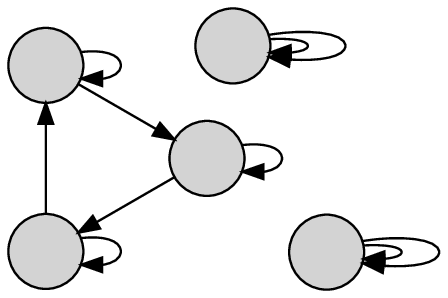}}

\fbox{\includegraphics[scale=0.3]{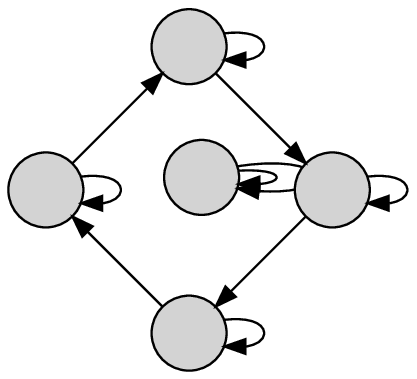}}
\fbox{\includegraphics[scale=0.3]{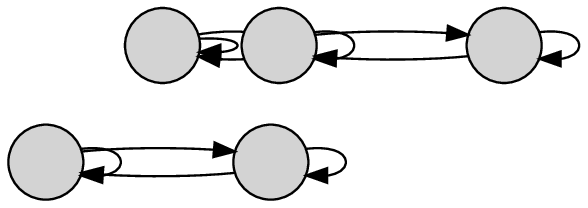}}
\fbox{\includegraphics[scale=0.3]{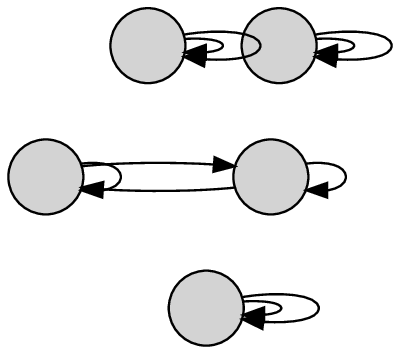}}
\fbox{\includegraphics[scale=0.3]{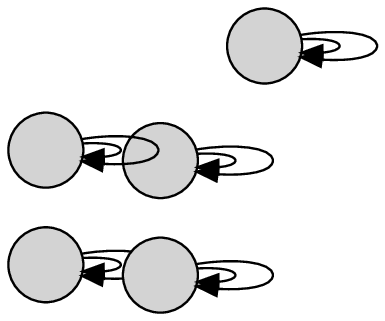}}

\end{widetext}

\section{Bijection with Lovelock Terms}

Bogdanos has pointed at a graphical representation of the terms of Lovelock's Lagrange density
\cite{LovelockJMP12,LovelockJMP13,BogdanosPRD70}:
A graph of a particular term which is a product of $n$ Riemann Tensors is constructed as follows:
(i) Start with a bottom row of $n$ nodes plus a top row of another $n$ nodes,
where nodes are associated left-to-right with a left-to-right reading of the $R$-factors of the product.
(ii) Add an edge from the $b$-th bottom node to the $t$-th top node if  a covariant lower index
of factor number $b$ appears as a contravariant index of factor number $t$.
This assignment works because the two contravariant indices are permutations of the covariant indices.
Since the Riemann Tensor is (skew)-symmetric in the first and in the last two indices, the
representation does not need to track the two edges at each node to distinguish a ``first'' from a ``second.''
Bogdanos' representation is a labeled, undirected 2-regular graph on $2n$ nodes.

Our graphical representation transforms these bipartite graphs once more \cite{CooperRSA4}:
\begin{itemize}
\item
Bogdanos' edges are turned into directed edges (arcs), always heading from a node of the bottom row to a
node on the top row.  
\item
Each pair of nodes associated with the same $R$-factor is 
collapsed
into
a single node, keeping all arcs fastened to their nodes. Loops appear if the $R$-factors
were already contracted (Ricci tensors).
\item
Labels are erased, meaning that the left-to-right reading orders of the product of the $R$
are all equivalent reflecting the usual commutative law for multiplications. 
The multiplicity may be recovered by examining the Automorphism Group of the new graph
\end{itemize}
This transformation of Bogdanos' bipartite graphs on $2n$ nodes to our 2-regular digraphs on $n$ nodes
is lossless (reversible).

\section{Line Graphs} 
Each digraph has an associated line digraph which is defined as follows \cite{HararyRendic9,SysloIPL15}:
Each arc in the digraph is rendered as one node in the line digraph.
An arc from node $a$ in the line digraph (corresponding to an arc from node $u$ to node $v$ in the digraph)
to node $b$ in the line digraph (corresponding to an arc from node $w$ to node $x$ in the digraph)
exists if the nodes $v$ and $w$ are the same, which means, if the two arcs in the graph
are adjacent to each other and have aligned directions.

Loops in the digraph create loops in the line digraph. A $k$-regular digraph on $n$ nodes
preduces a $k$-regular line digraph on $kn$ nodes.

\section{Single Edge Cuts}\label{sec.scut}
If a single arc of the 2-regular digraphs discussed above is cut,
a pair of two additional nodes appears at that place: one with indegree 1 and outdegree 0,
another with indegree 0 and outdegree 1; the other nodes keep their degrees. 
The equivalent action in the tensor product 
is to leave one pair of indices uncontracted; the equivalent action in a $\varphi^4$-theory
is to move from the vacuum diagrams to the diagrams of Green's functions.
In a preparational step one might have enumerated the marked line graphs where
the marked nodes in the line graphs would  indicate which lines are cut in
their base graphs. 

The notation $U_{k,+2}(n,c)$ will be used for the number the unlabeled digraphs
with $c$ components 
with $n-2$ nodes with indegree and outdegree $k$ plus 2 nodes with the mixed in-outdegrees
of (0,1) and (1,0) created by the cut.
(Adopting the nomenclature of  Chae et al. \cite{ChaeDM307}, these ``almost'' regular graphs will
be called fairly regular graphs.)
The enumeration for any number of components
is $U_{k,+2}(n)\equiv \sum_c U_{k,+1}(n,c)$.
Only the connected graphs (where $c=1$) are relevant because the graphs with 
more than one component can be constructed by taking one connected graph 
and
adding any number of graphs of the 
$k$-regular class. So the ordinary generating function
of the fairly regular digraphs $U_{k,+2}(n)$ 
is the product of the ordinary generating function
of $U_{k,+2}(n,1)$ by the ordinary generating function of $U_{k}(n)$,
and the equivalent statement holds for the
exponential generating function of the labeled graphs.

The full set of these 1-in 1-out graphs up to 6 nodes obtained
by cutting an arc follows.
The cases without multiarcs up to 7 nodes I have already shown in
a manuscript in \cite[A005642]{sloane}.

\begin{widetext}

\subsection{1 connected graph on 2 nodes}
\includegraphics[scale=0.3]{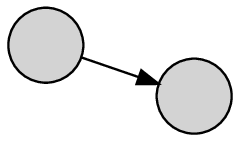}
\subsection{1 connected graph on 3 nodes}
\includegraphics[scale=0.3]{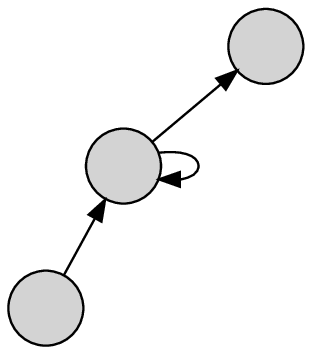}
\subsection{3 connected graphs on 4 nodes}
\includegraphics[scale=0.3]{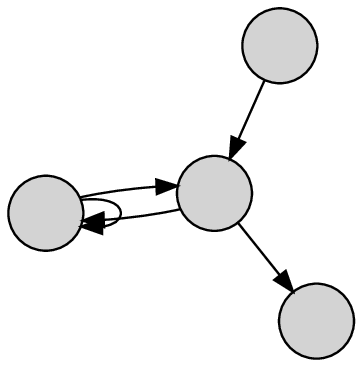}
\includegraphics[scale=0.3]{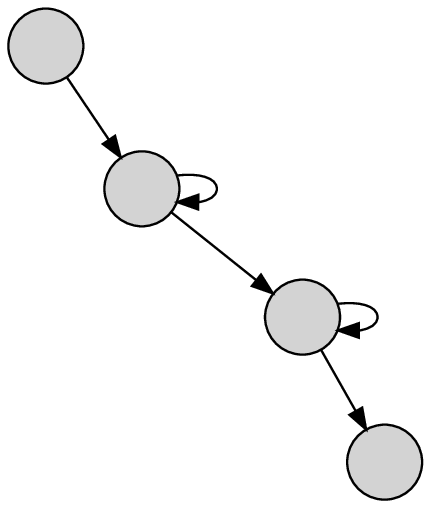}
\includegraphics[scale=0.3]{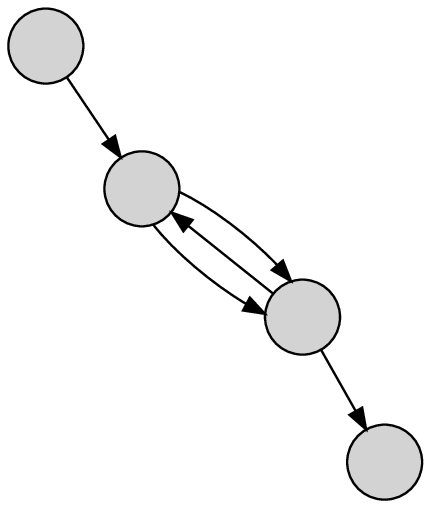}
\subsection{12 connected graphs on 5 nodes}
\includegraphics[scale=0.3]{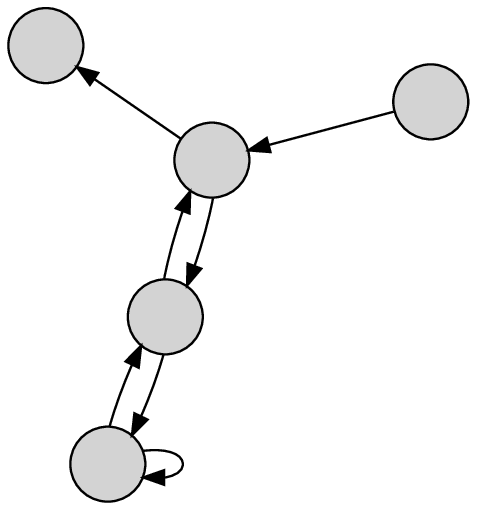}
\includegraphics[scale=0.3]{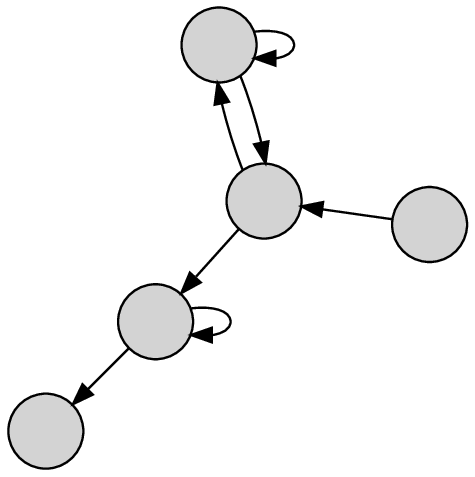}
\includegraphics[scale=0.3]{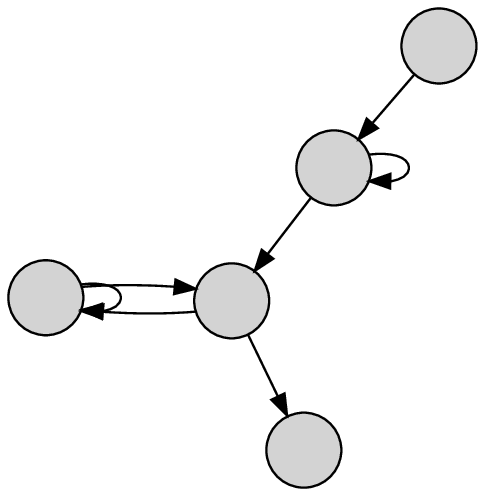}
\includegraphics[scale=0.3]{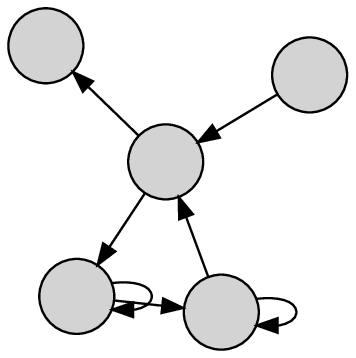}
\includegraphics[scale=0.3]{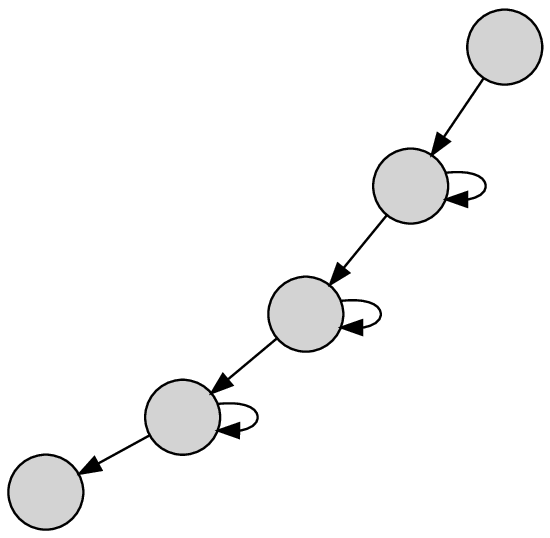}
\includegraphics[scale=0.3]{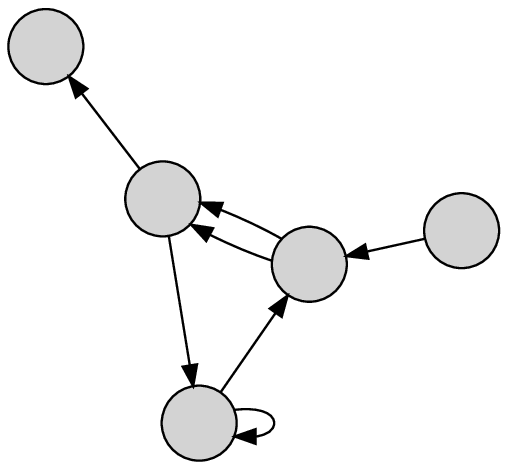}
\includegraphics[scale=0.3]{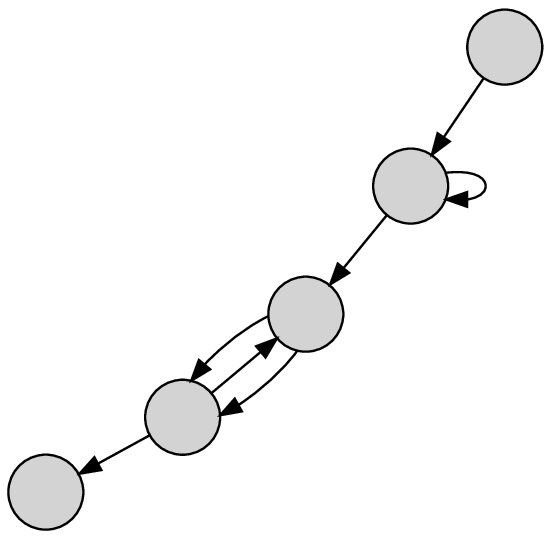}
\includegraphics[scale=0.3]{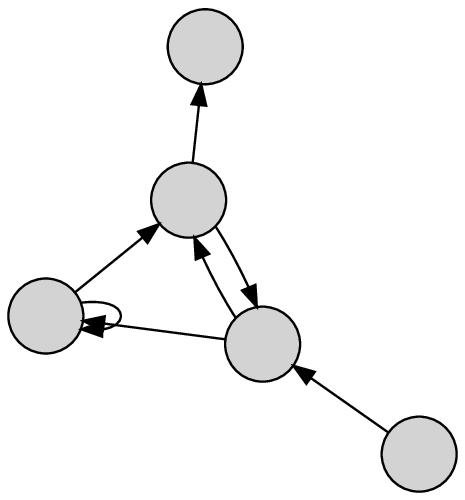}
\includegraphics[scale=0.3]{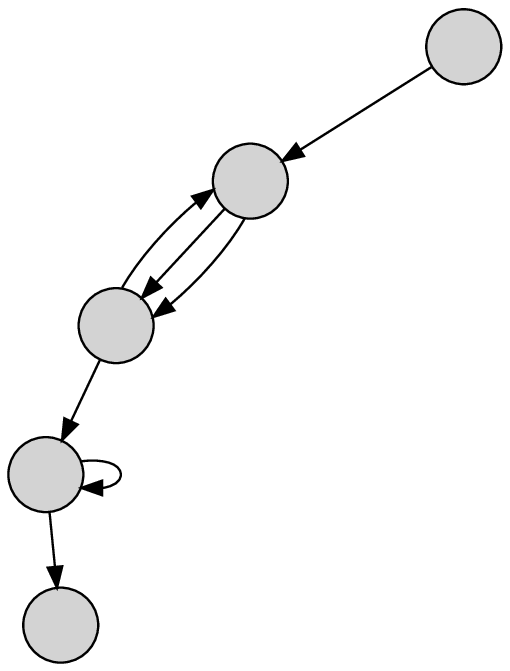}
\includegraphics[scale=0.3]{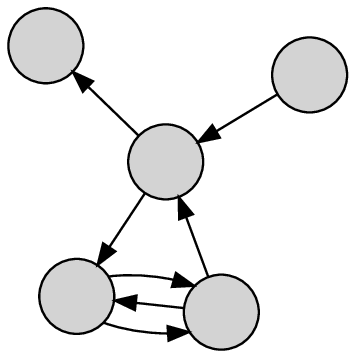}
\includegraphics[scale=0.3]{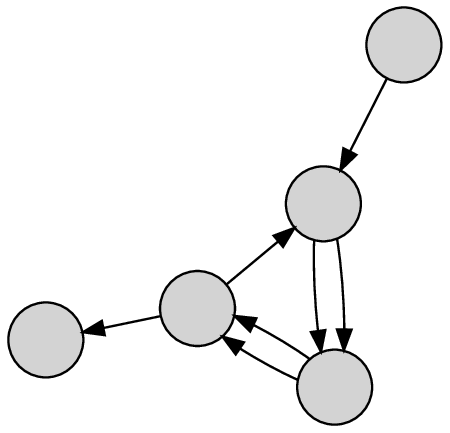}
\includegraphics[scale=0.3]{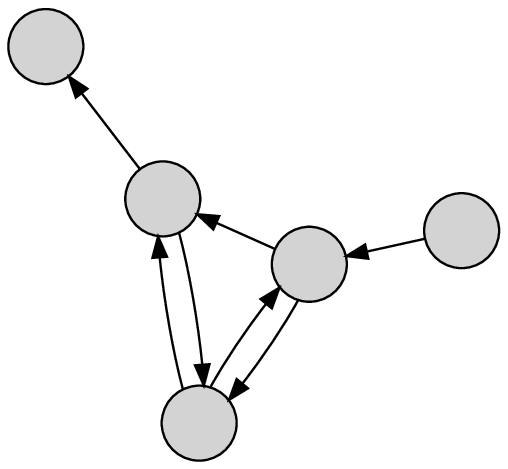}
\subsection{61 connected graphs on 6 nodes}
\includegraphics[scale=0.3]{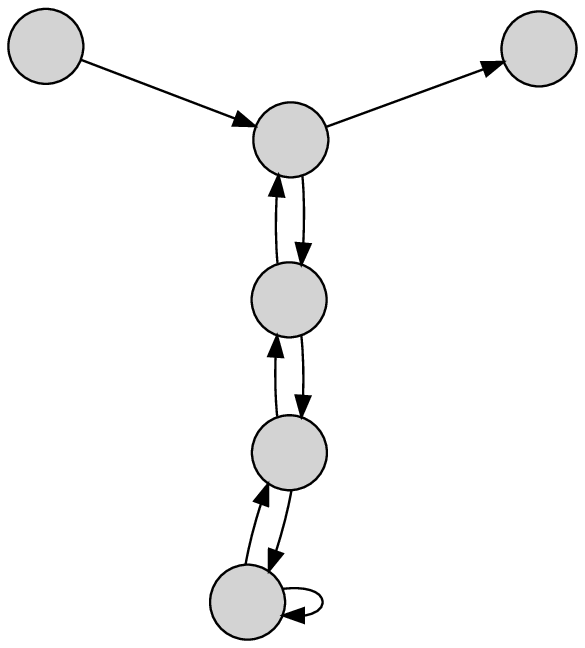}
\includegraphics[scale=0.3]{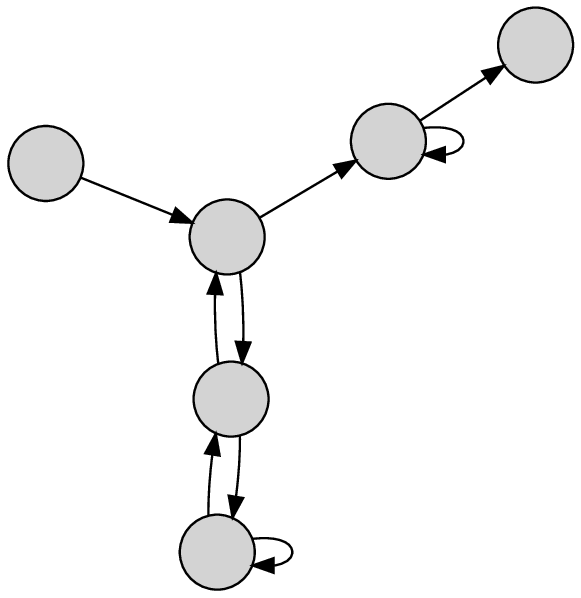}
\includegraphics[scale=0.3]{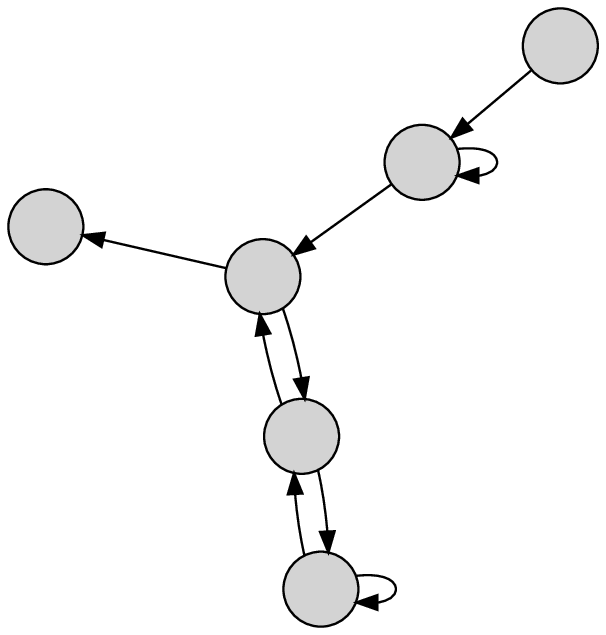}
\includegraphics[scale=0.3]{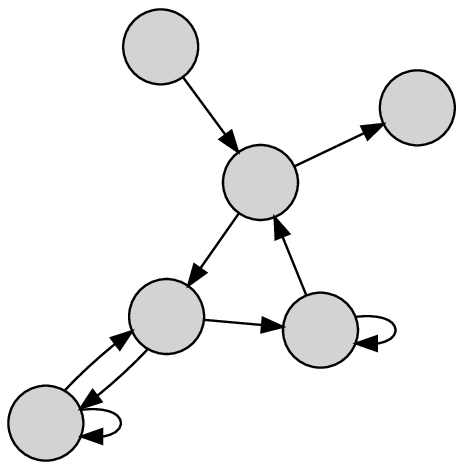}
\includegraphics[scale=0.3]{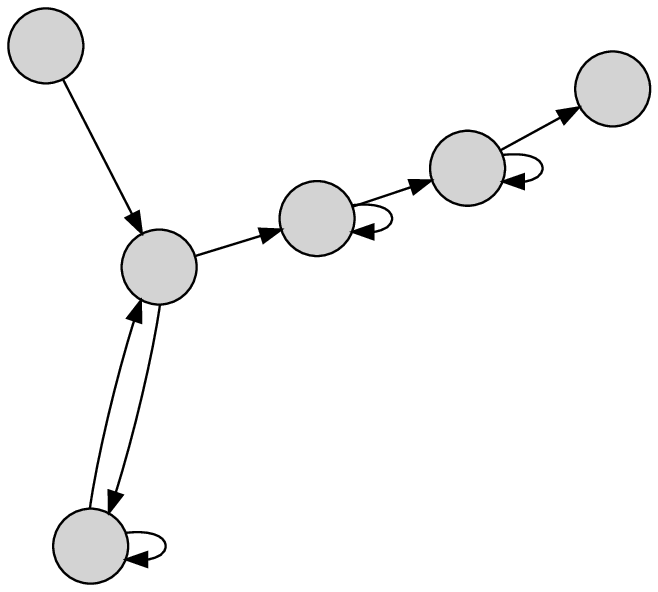}
\includegraphics[scale=0.3]{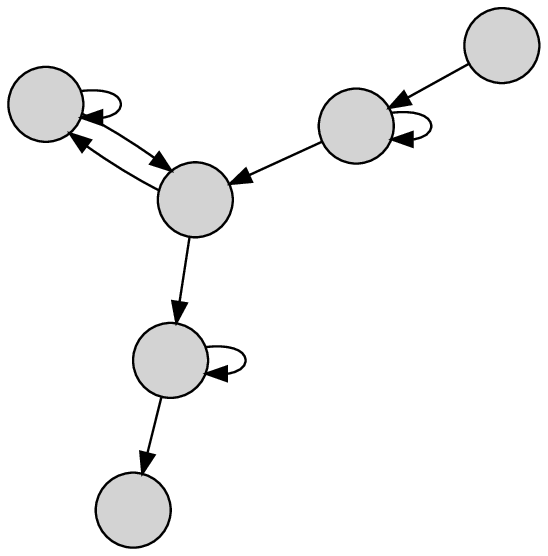}
\includegraphics[scale=0.3]{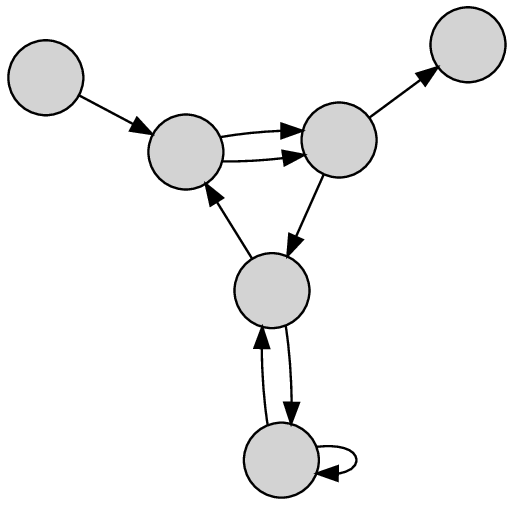}
\includegraphics[scale=0.3]{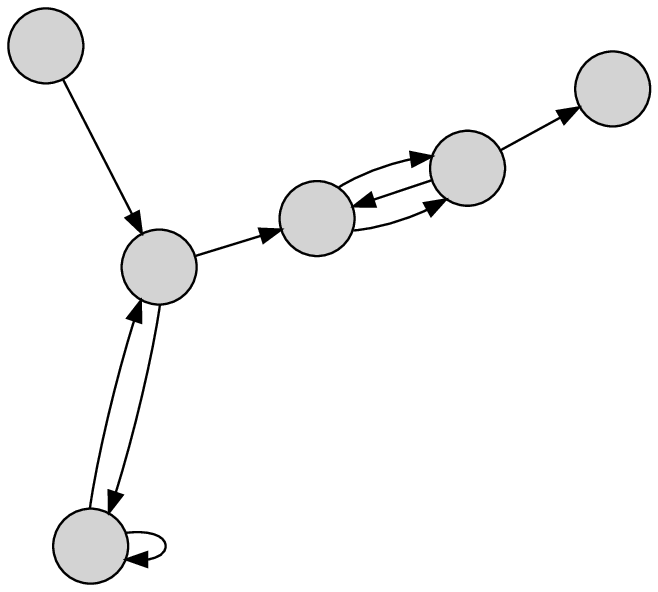}
\includegraphics[scale=0.3]{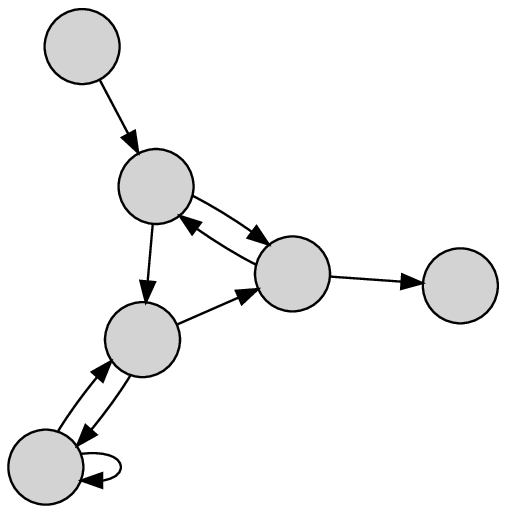}
\includegraphics[scale=0.3]{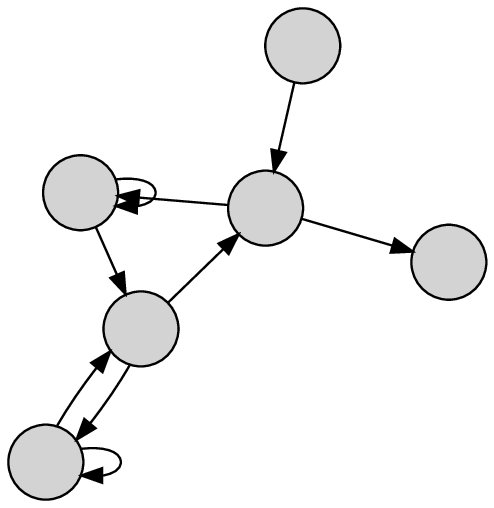}
\includegraphics[scale=0.3]{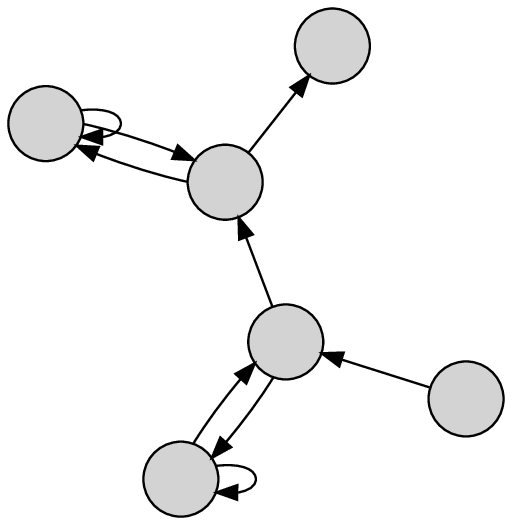}
\includegraphics[scale=0.3]{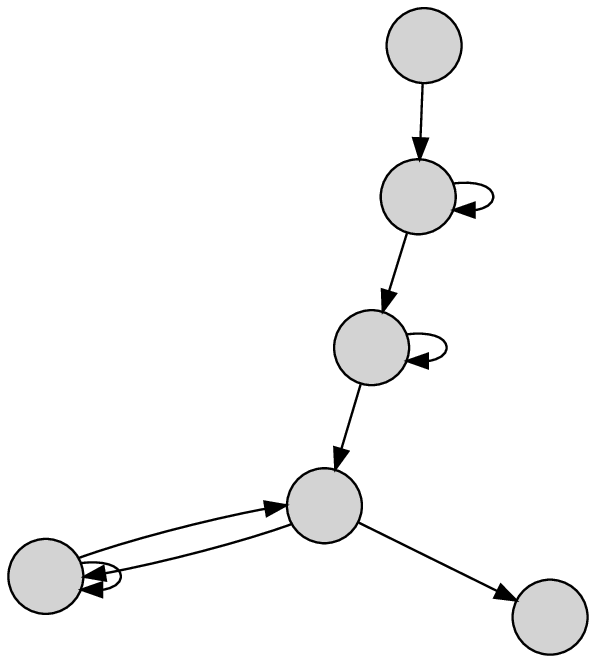}
\includegraphics[scale=0.3]{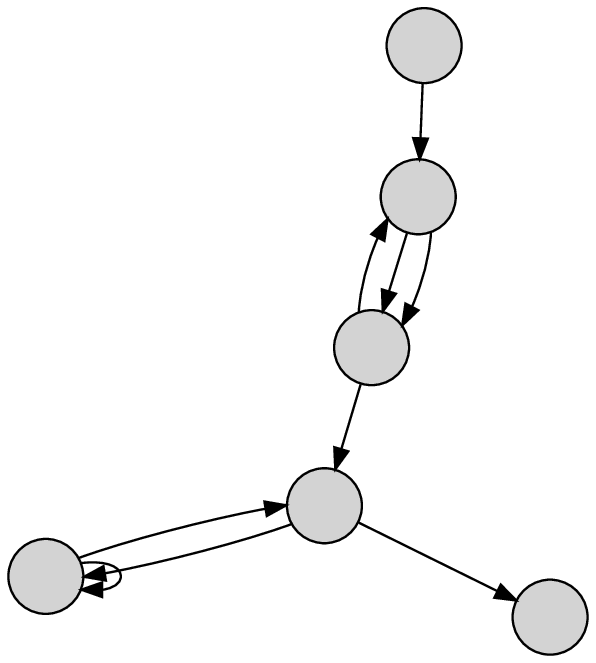}
\includegraphics[scale=0.3]{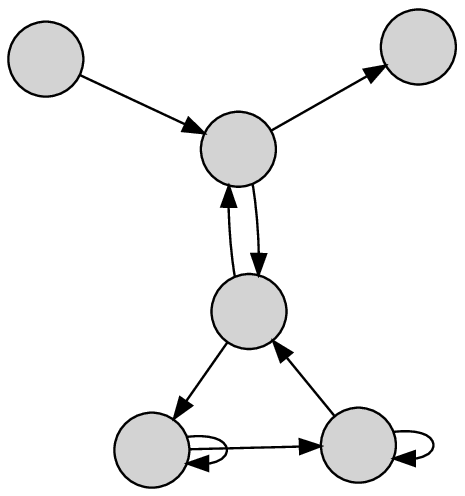}
\includegraphics[scale=0.3]{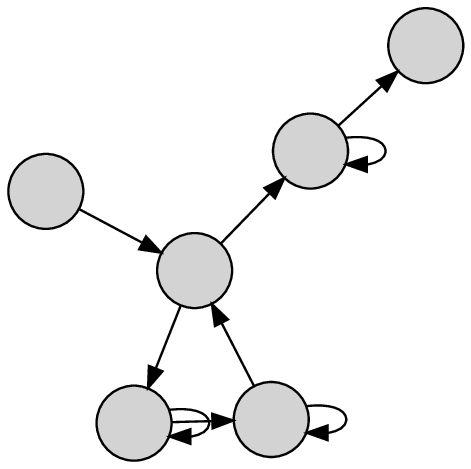}
\includegraphics[scale=0.3]{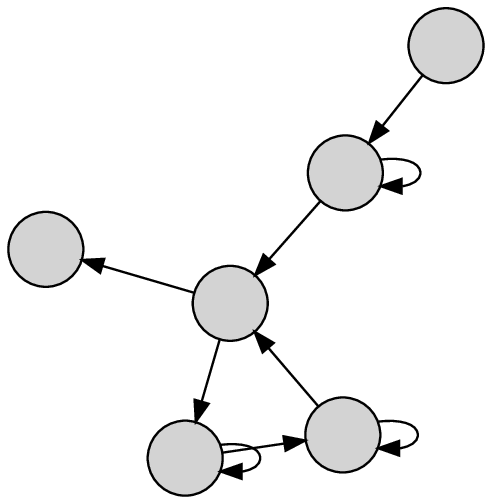}
\includegraphics[scale=0.3]{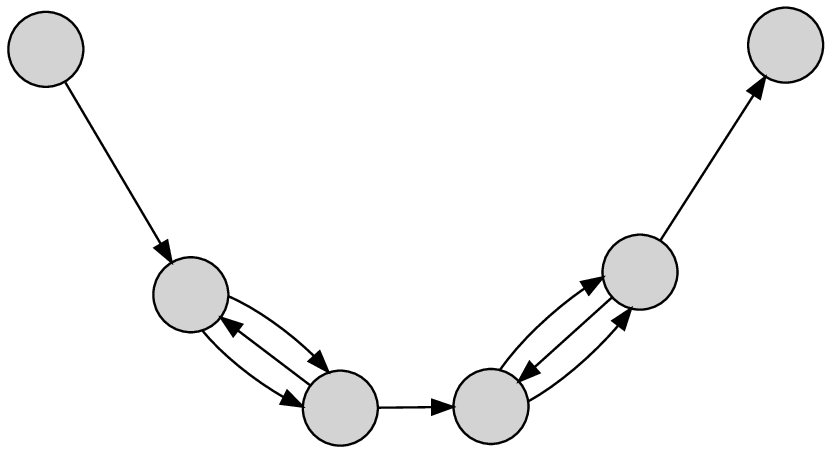}
\includegraphics[scale=0.3]{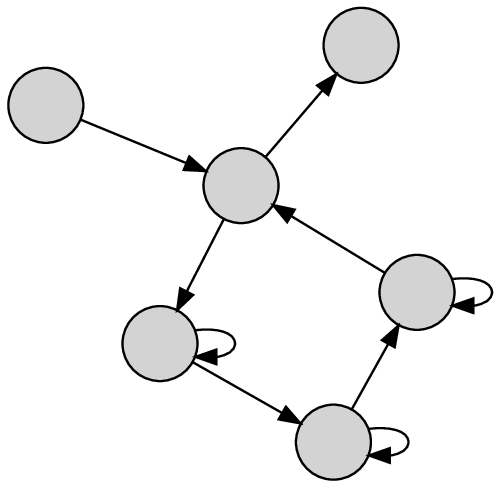}
\includegraphics[scale=0.3]{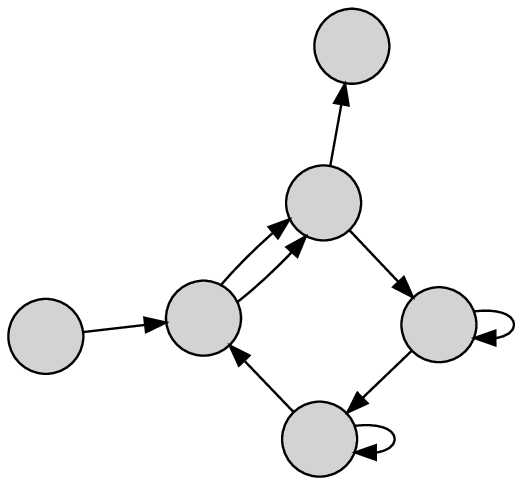}
\includegraphics[scale=0.3]{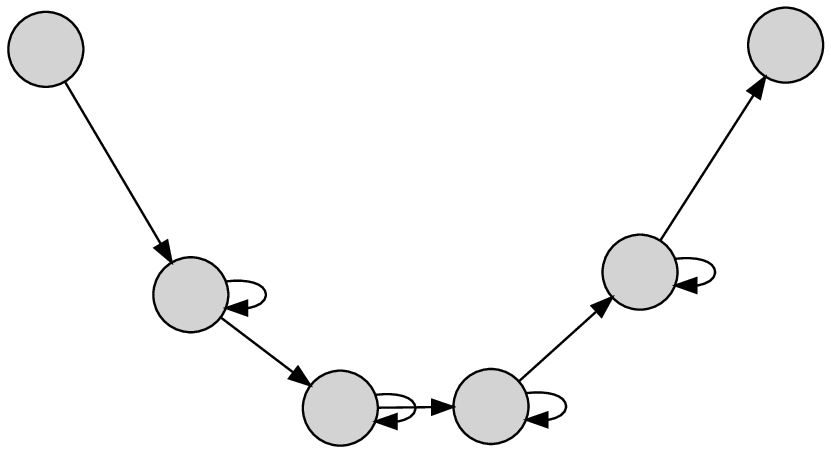}
\includegraphics[scale=0.3]{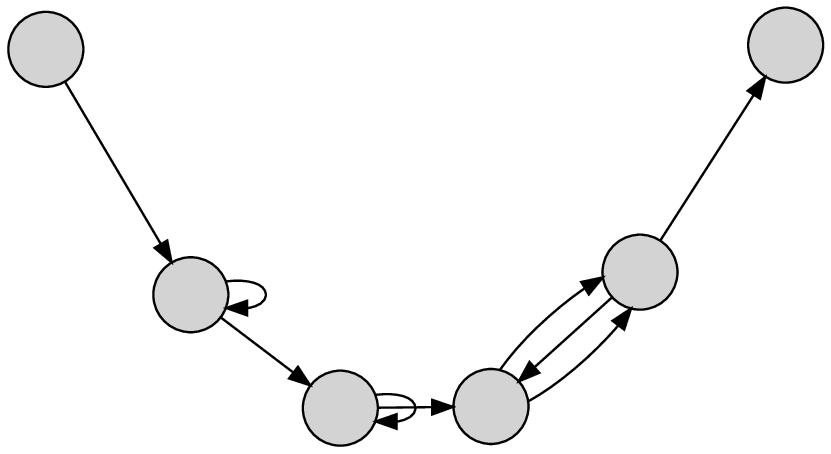}
\includegraphics[scale=0.3]{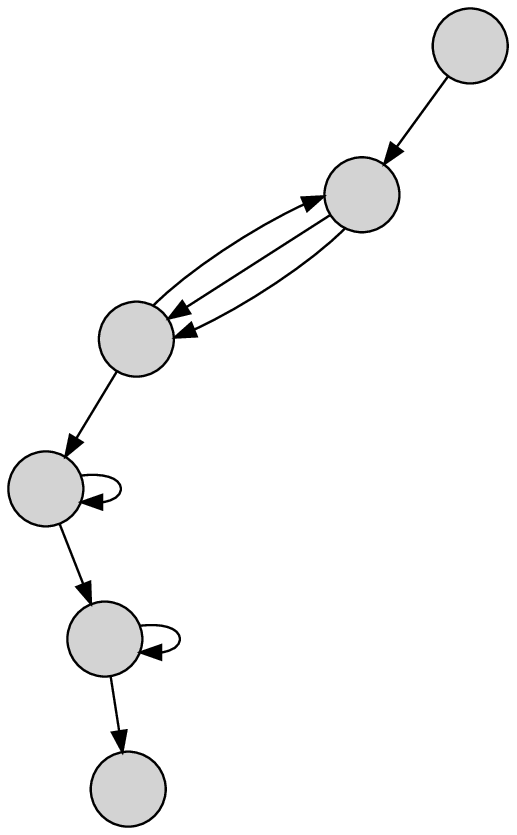}
\includegraphics[scale=0.3]{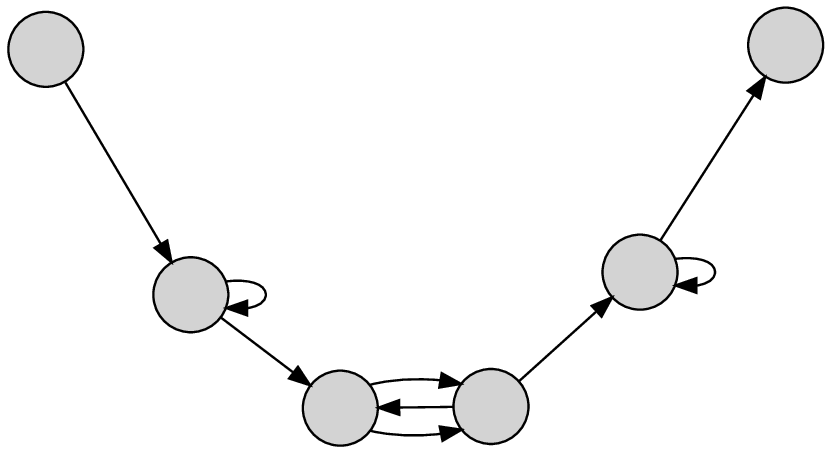}
\includegraphics[scale=0.3]{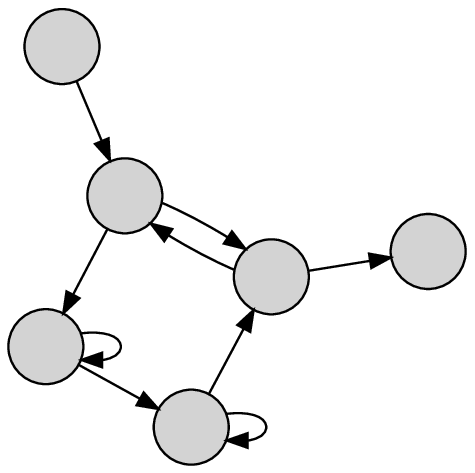}
\includegraphics[scale=0.3]{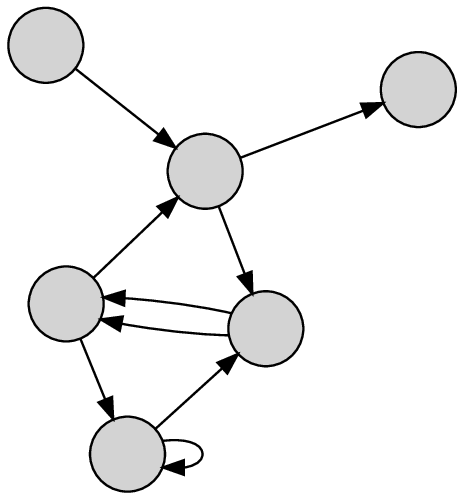}
\includegraphics[scale=0.3]{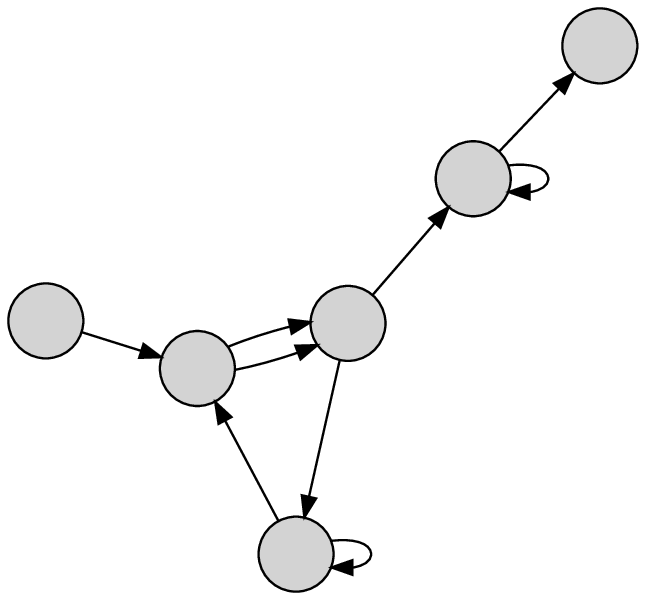}
\includegraphics[scale=0.3]{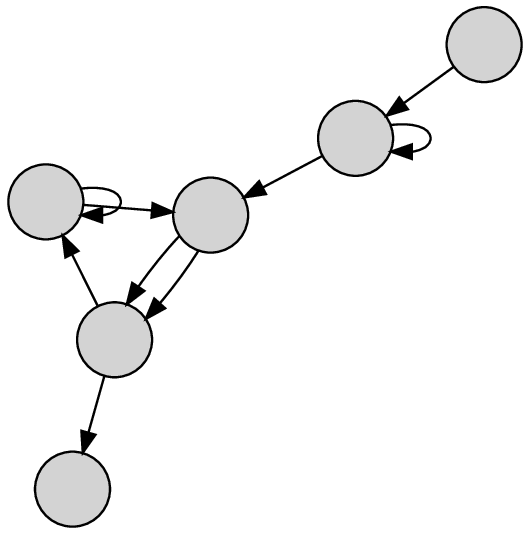}
\includegraphics[scale=0.3]{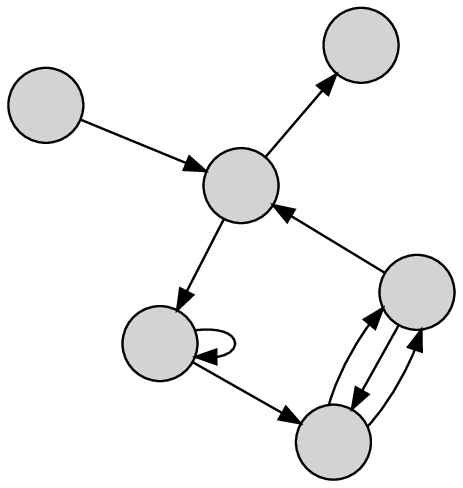}
\includegraphics[scale=0.3]{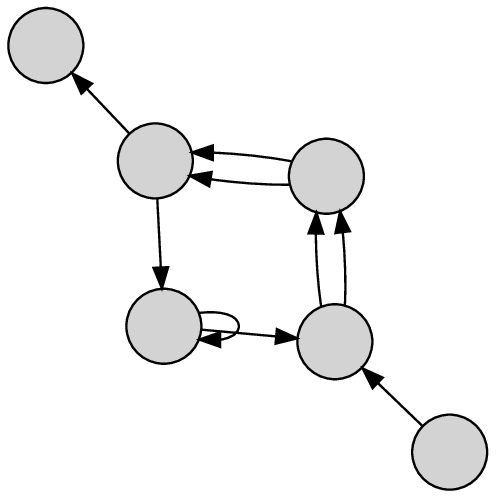}
\includegraphics[scale=0.3]{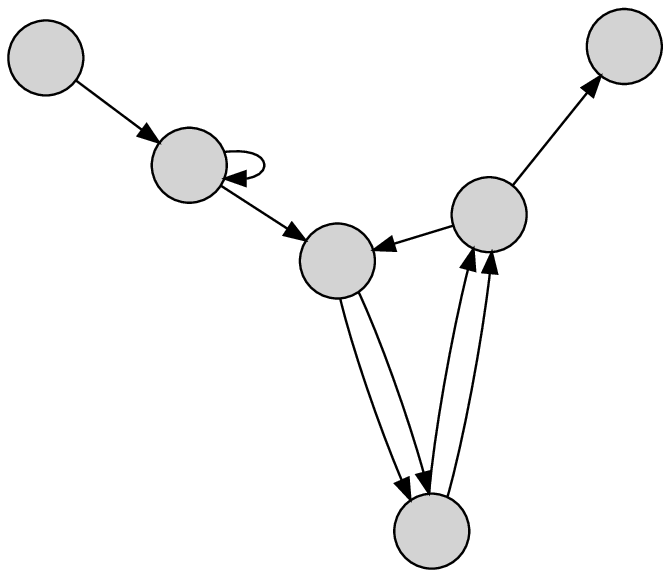}
\includegraphics[scale=0.3]{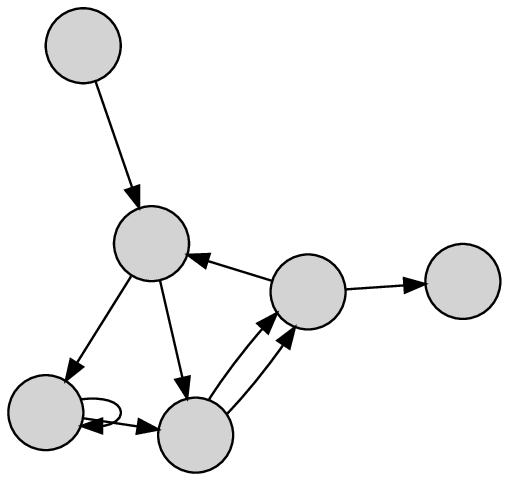}
\includegraphics[scale=0.3]{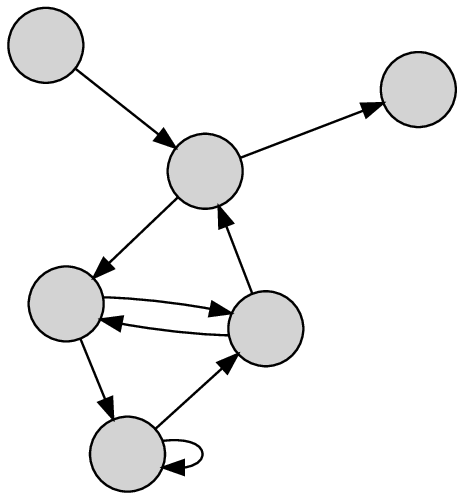}
\includegraphics[scale=0.3]{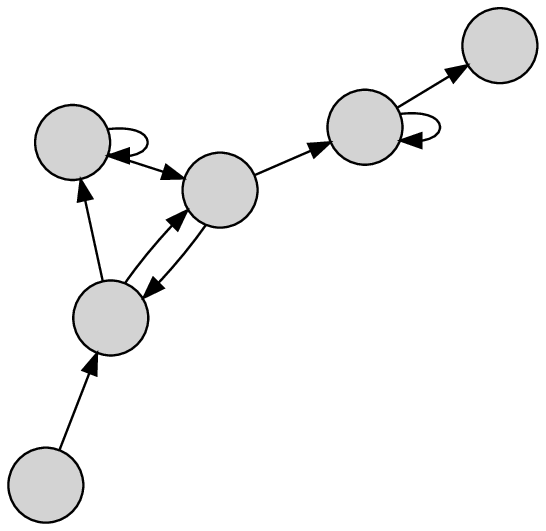}
\includegraphics[scale=0.3]{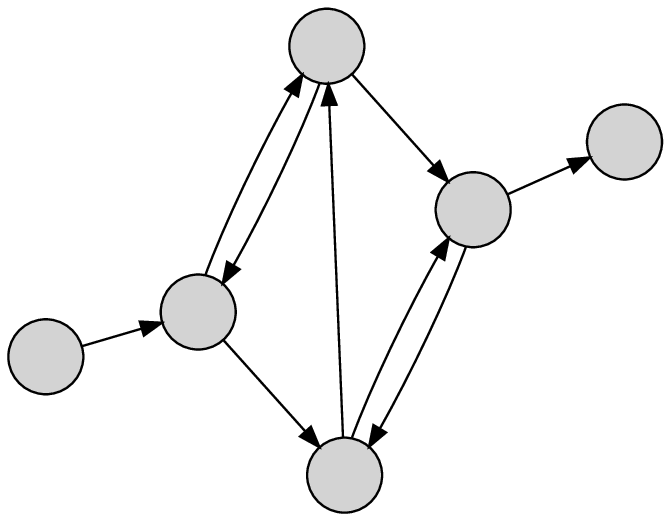}
\includegraphics[scale=0.3]{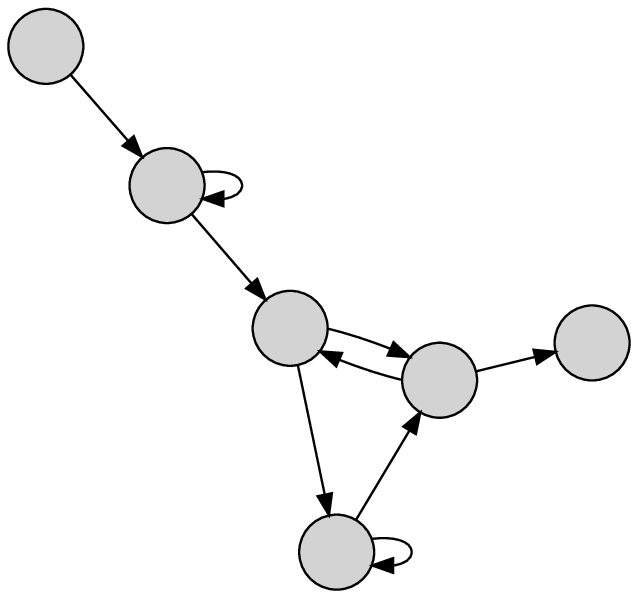}
\includegraphics[scale=0.3]{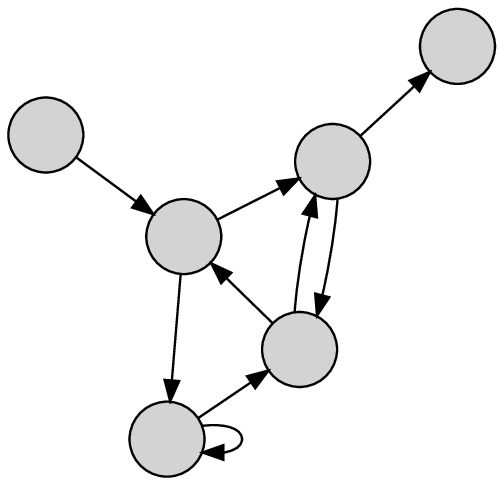}
\includegraphics[scale=0.3]{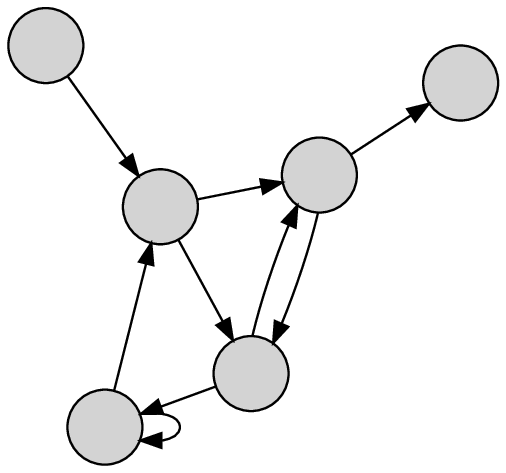}
\includegraphics[scale=0.3]{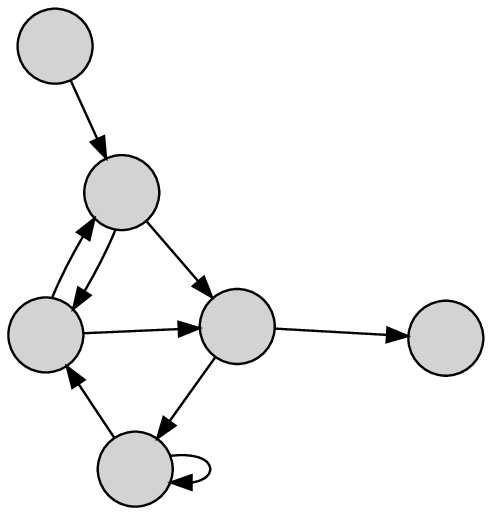}
\includegraphics[scale=0.3]{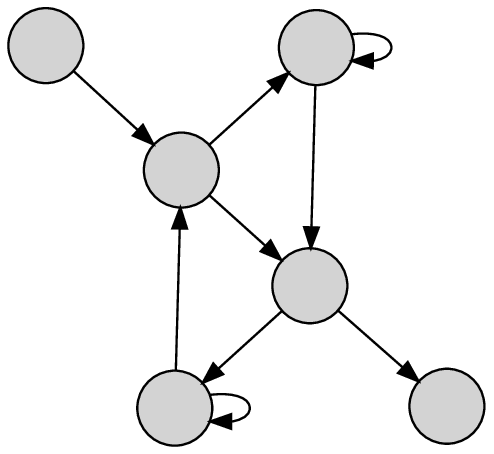}
\includegraphics[scale=0.3]{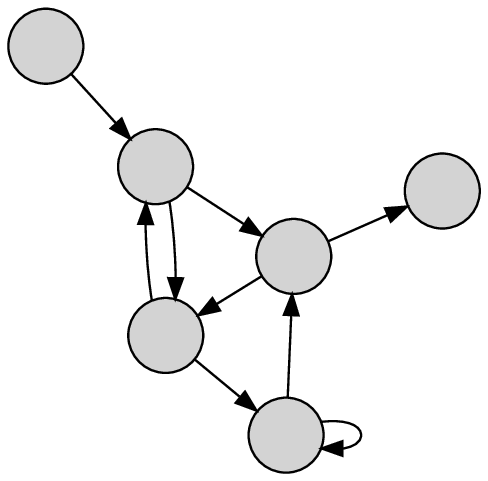}
\includegraphics[scale=0.3]{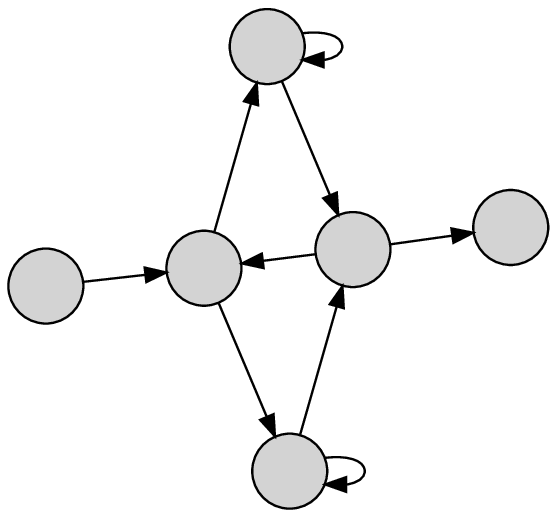}
\includegraphics[scale=0.3]{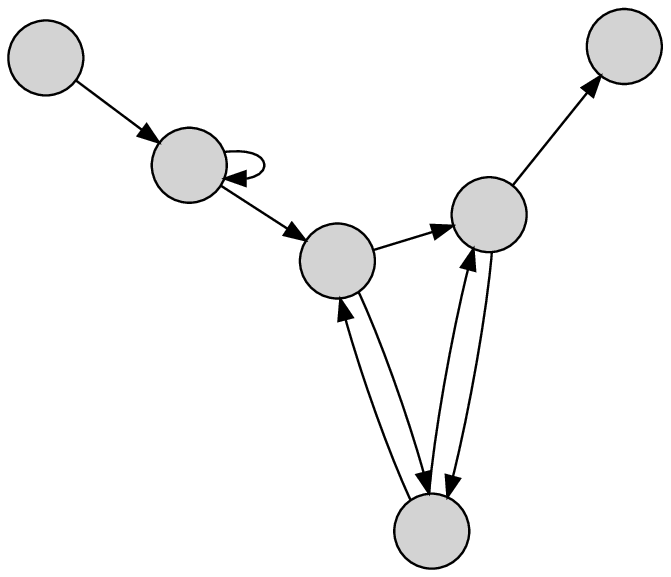}
\includegraphics[scale=0.3]{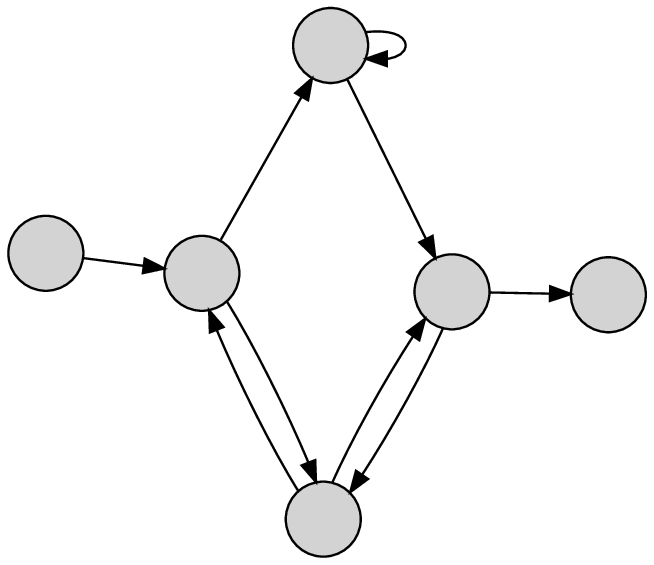}
\includegraphics[scale=0.3]{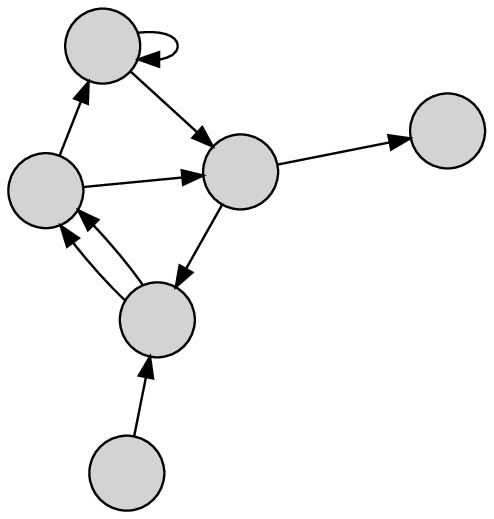}
\includegraphics[scale=0.3]{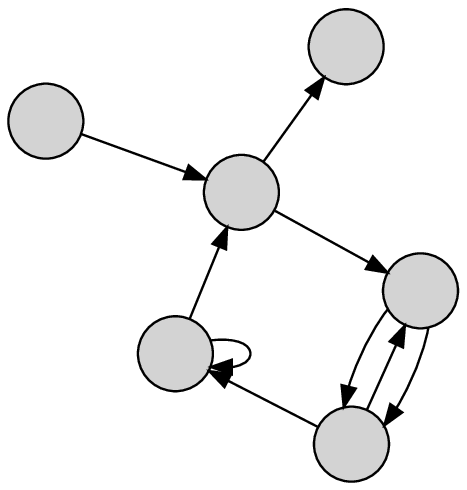}
\includegraphics[scale=0.3]{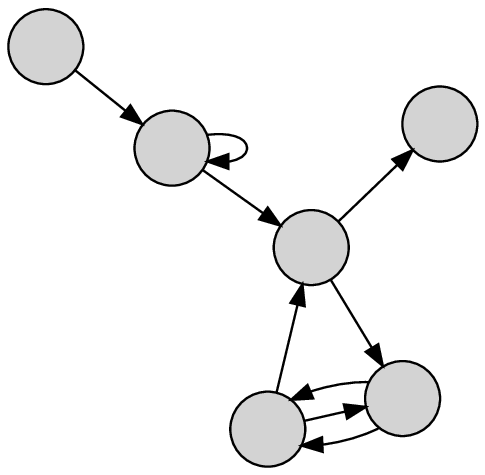}
\includegraphics[scale=0.3]{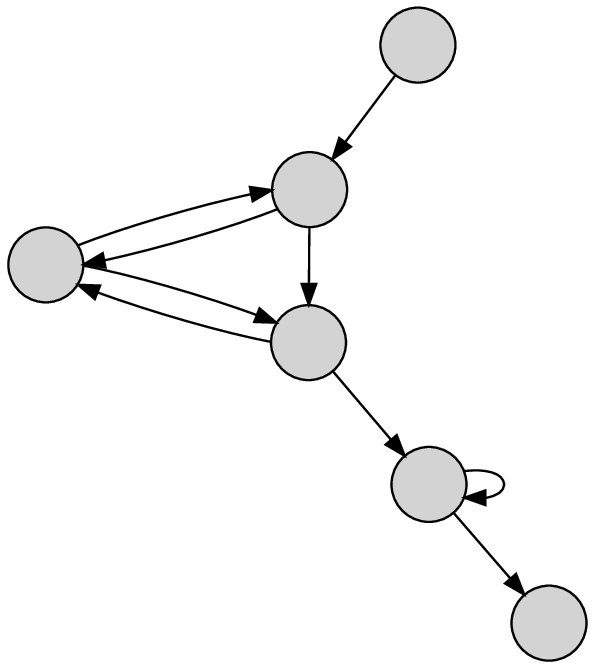}
\includegraphics[scale=0.3]{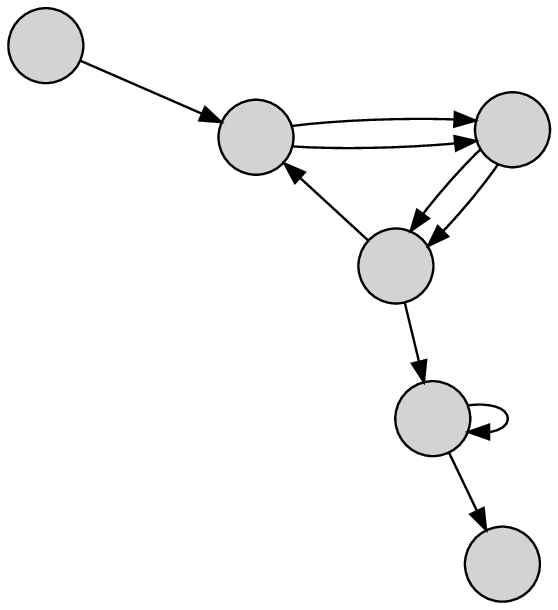}
\includegraphics[scale=0.3]{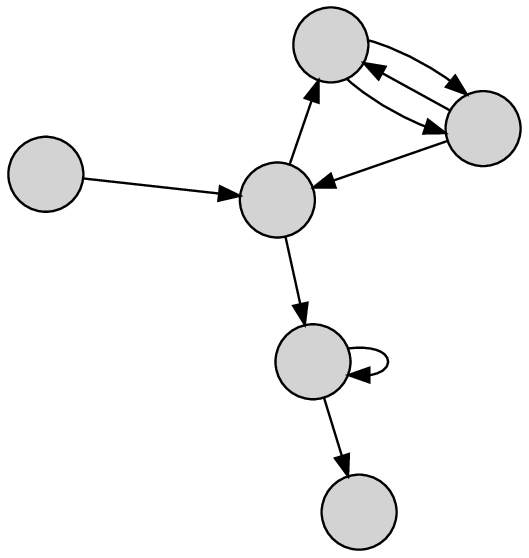}
\includegraphics[scale=0.3]{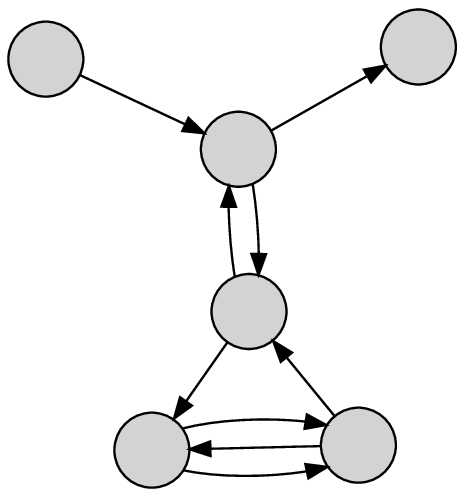}
\includegraphics[scale=0.3]{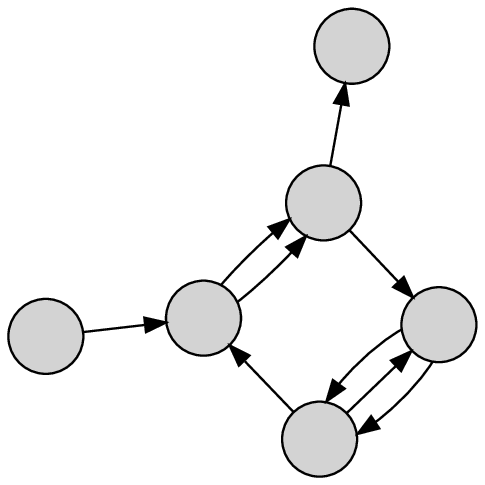}
\includegraphics[scale=0.3]{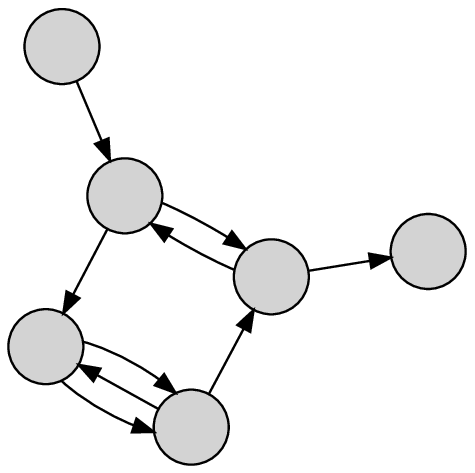}
\includegraphics[scale=0.3]{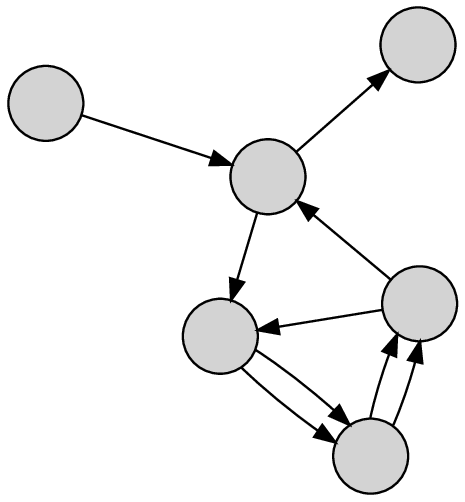}
\includegraphics[scale=0.3]{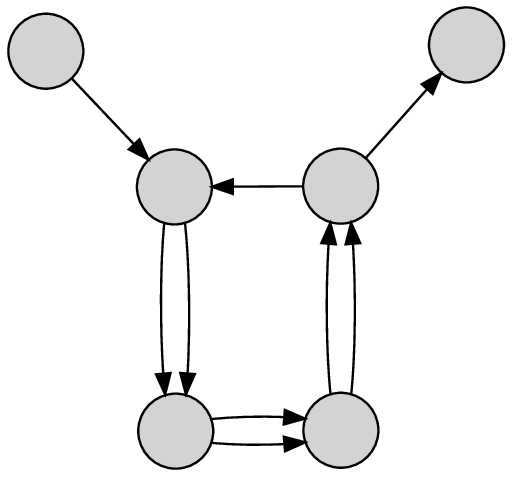}
\includegraphics[scale=0.3]{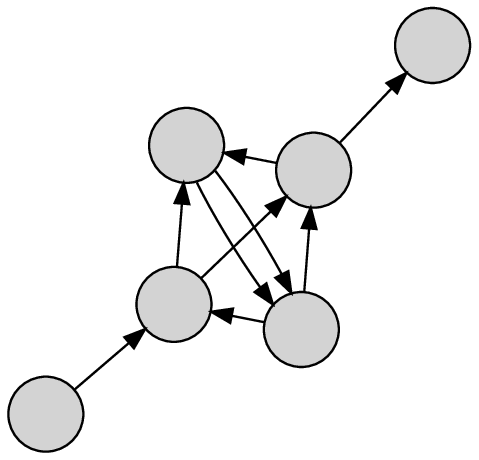}
\includegraphics[scale=0.3]{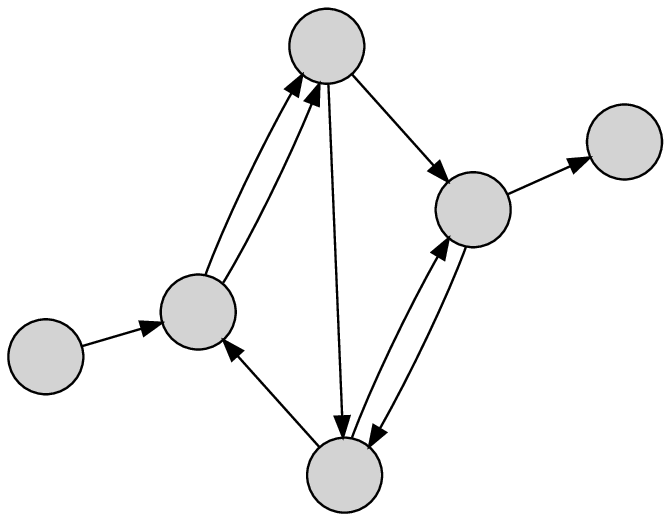}
\includegraphics[scale=0.3]{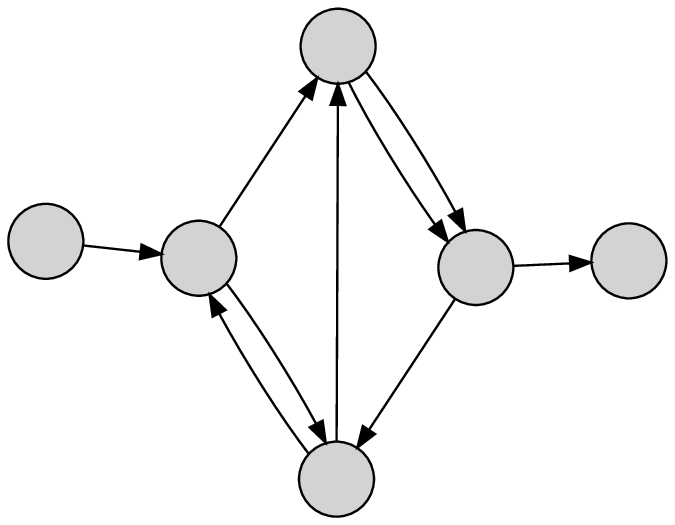}
\includegraphics[scale=0.3]{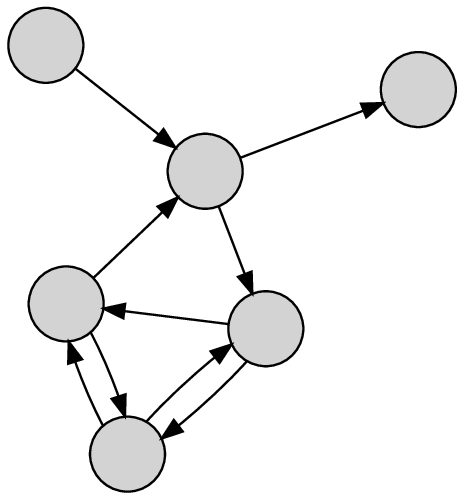}
\includegraphics[scale=0.3]{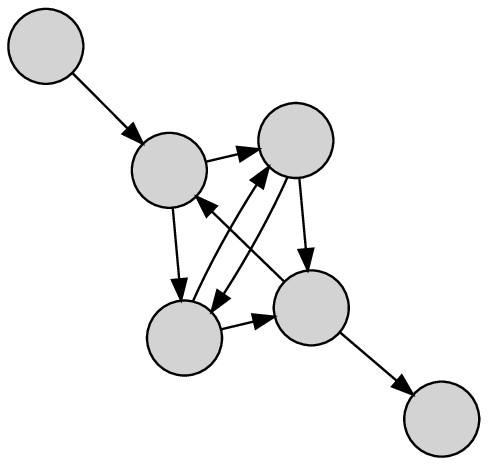}
\includegraphics[scale=0.3]{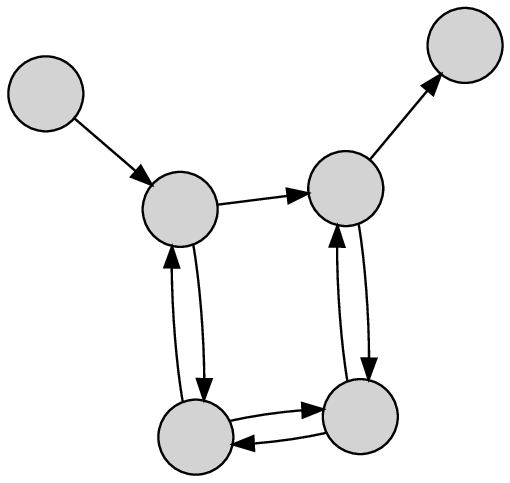}
\includegraphics[scale=0.3]{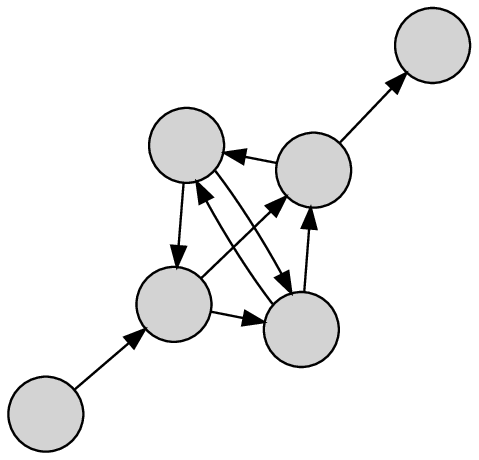}

\end{widetext}

A summary of the enumeration with one arc cut
is in Table \ref{tab.f}.
\begin{table}[htb]
\begin{ruledtabular}
\begin{tabular}{r|rrr}
$n$ & $U_{2,+2}(n,1)$ & $L_{2,+2}(n,1)$ &$L_{2,+2}(n)$ \\
\hline
2 & 1 & 2& 2 \\
3 & 1& 6& 12  \\
4 & 3 & 72& 132  \\
5 & 12 & 1440 & 2400  \\
6 & 61 & 43200& 66060  \\
7 & 365 & 1806840& 2572920 \\
\end{tabular}
\end{ruledtabular}
\caption{Connected unlabeled 2-regular digraphs after cutting one arc, and 
corresponding labeled 2-regular digraphs (one or any number of components)
allowing loops and multiarcs.} 
\label{tab.f}
\end{table}

\section{Pairs of Cuts}\label{sec.2cut}
Cutting two distinct arcs of a connected 2-regular digraph
produces $U_{2,+4}(n,c)$ unlabeled digraphs with $c$ components
where 4 nodes
appear at the cuts: 2 with indegree 1 and outdegree 0,
and another 2 with indegree 0 and outdegree 1. The $n-4$ ``internal''
nodes keep their degrees.
(Cutting the same arc twice is not interesting because that merely
produces two disconnected graphs already covered by the previous section.)
The operation with two cuts is topologically more complex
than products of a single cut. (i) There is no longer
a 1-to-1 correspondence from the graph with the 2-in 2-out nodes
to the uncut parent. There are generally two variants
of merging the nodes with total degree 1 to
rebuild a 2-regular digraph; exceptions are cases with high symmetry
like in Section \ref{sec.2c5} which obviously has the parent shown
in Section \ref{sec.1}.
(ii) If the two cuts target a pair of arcs between a single pair
of nodes or more generally two arcs of a 2-connected graph, 
a graph with more than one component (i.e. two
graphs of the type in Section \ref{sec.scut}) may result.
So the generating function of $\sum_c U_{2,+4}(n,c)$ is the product
of the generating functions $U_{2,+4}(n,1)$ and $\sum_c U_{2}(n,c)$
to account for graphs where all four nodes of degree 1 are in
the same component, \emph{plus} the symmetrized square
generating function of $\sum_c U_{2,2}(n,c)$ by $\sum_c U_2(n,c)$
to account for the
graphs where two nodes of degree 1 are in one component,
the other two nodes in another component, and any other
number of 2-regular graphs in further components.
[The symmetrized square of a generating function $g(x)$ is
$[g(x)^2+g(x^2)]/2$ following from the cycle index for the symmetric group $S_2$
\cite[I60]{Flajolet}.]

A gallery of these 2-in 2-out graphs up to 7 nodes illustrates
that concept in the next 3 subsections, showing only
the connected digraphs. 
\subsection{1 connected graph on 5 nodes} \label{sec.2c5}
\includegraphics[scale=0.25]{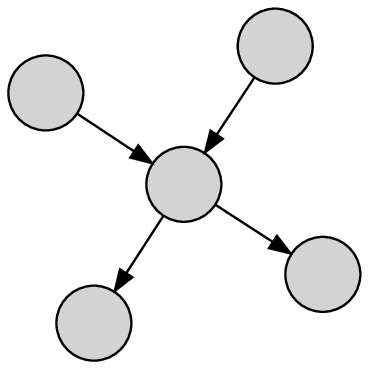}
\subsection{4 connected graphs on 6 nodes}
\includegraphics[scale=0.25]{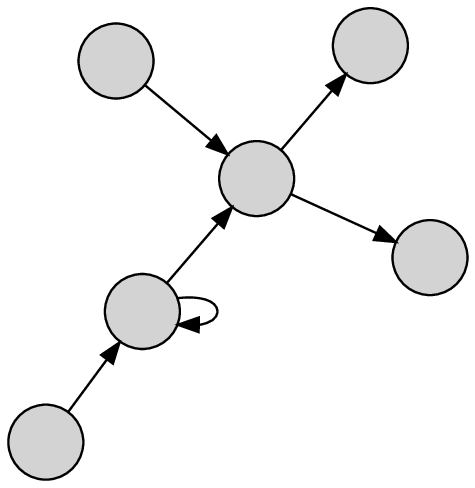}
\includegraphics[scale=0.25]{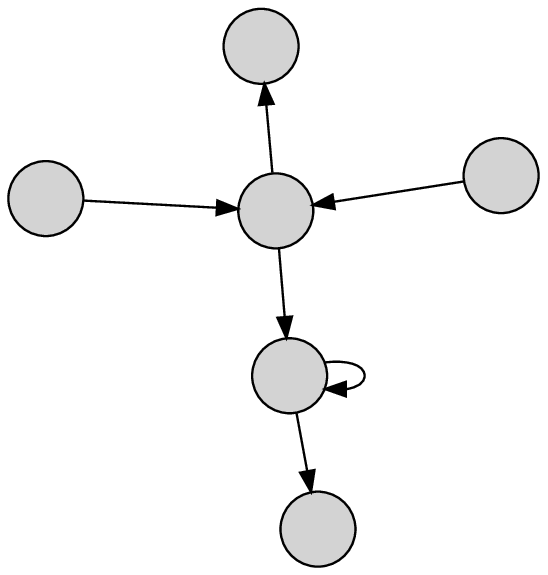}
\includegraphics[scale=0.25]{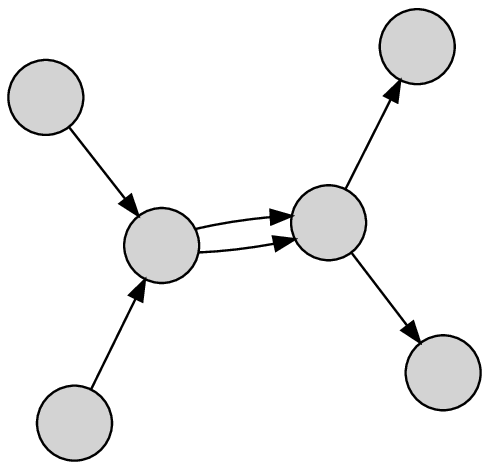}
\includegraphics[scale=0.25]{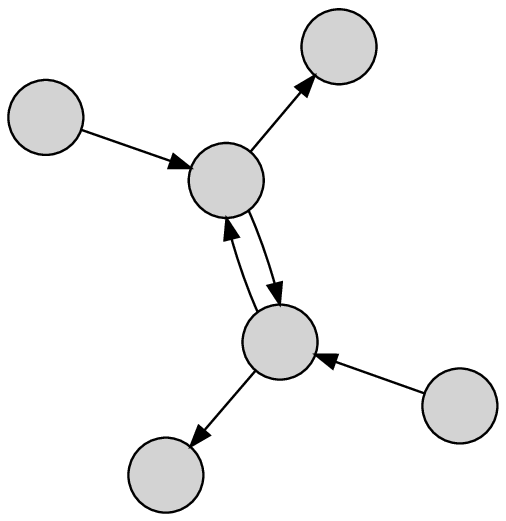}
\subsection{21 connected graphs on 7 nodes}
\includegraphics[scale=0.25]{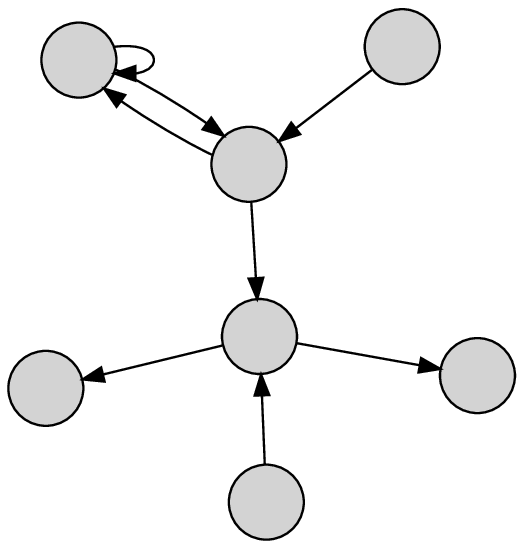}
\includegraphics[scale=0.25]{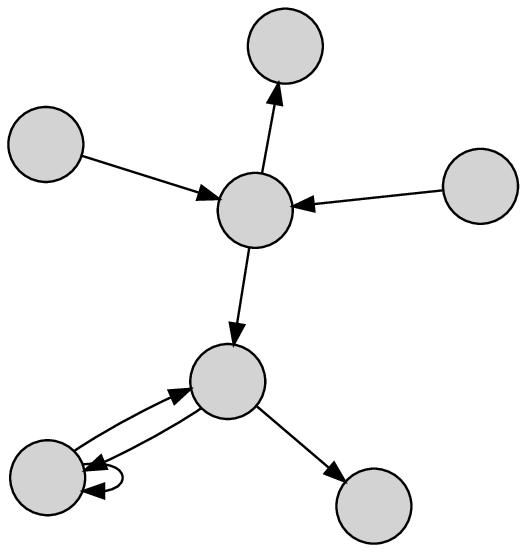}
\includegraphics[scale=0.25]{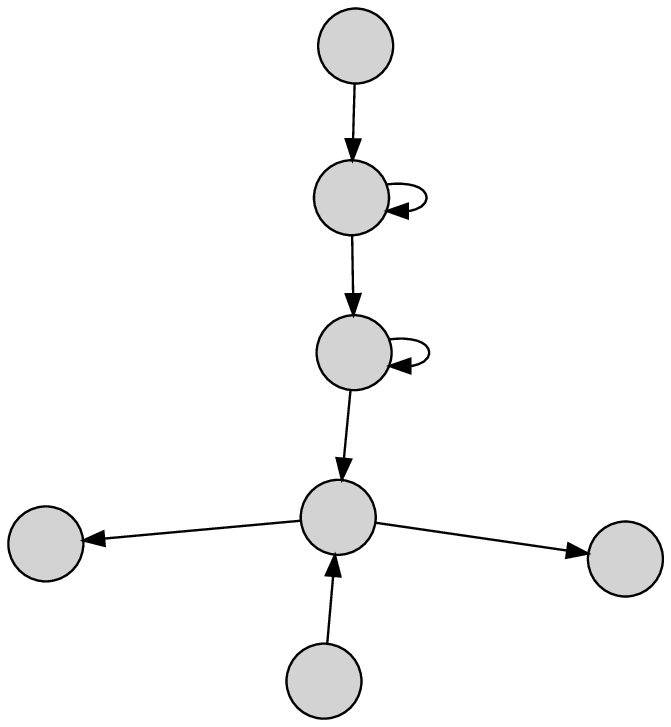}
\includegraphics[scale=0.25]{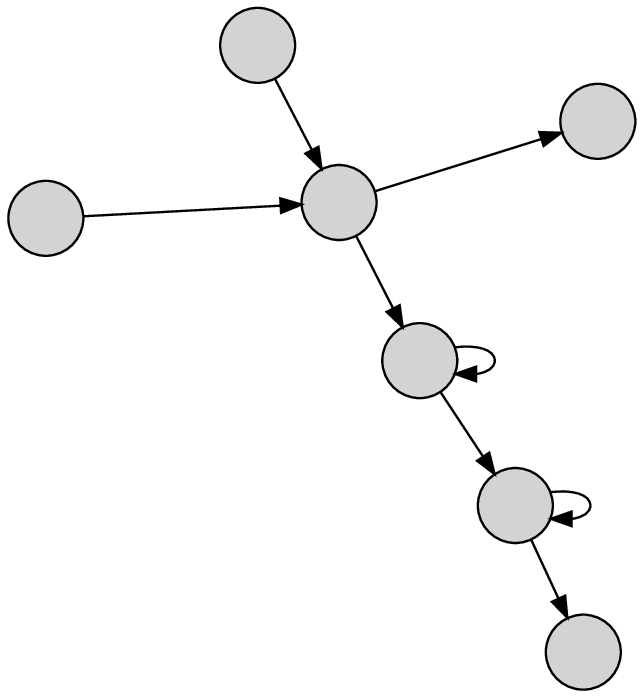}
\includegraphics[scale=0.25]{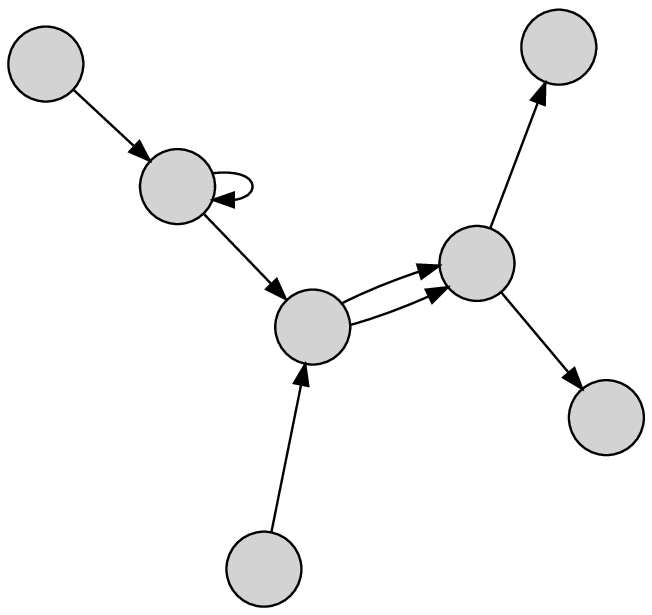}
\includegraphics[scale=0.25]{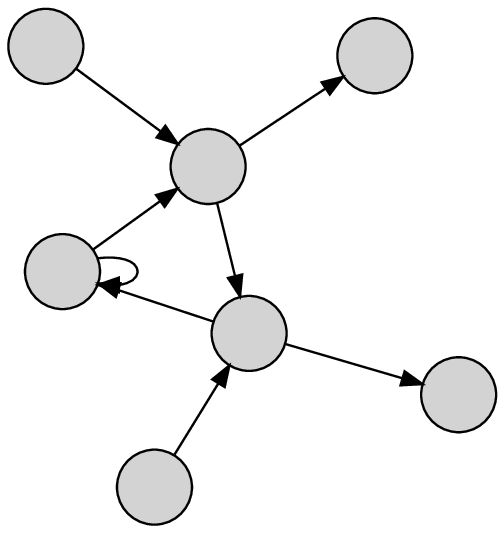}
\includegraphics[scale=0.25]{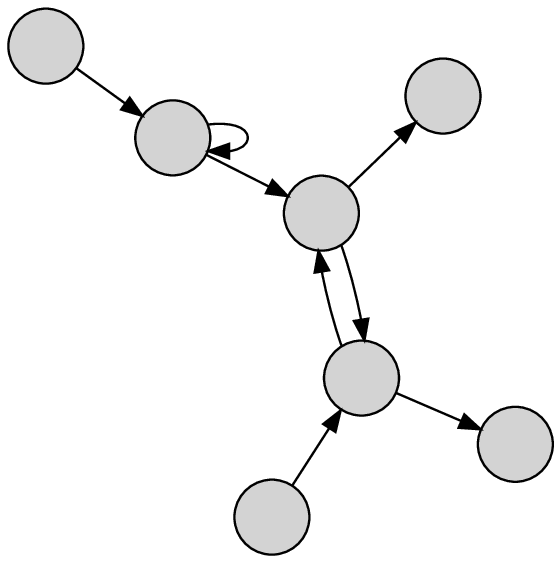}
\includegraphics[scale=0.25]{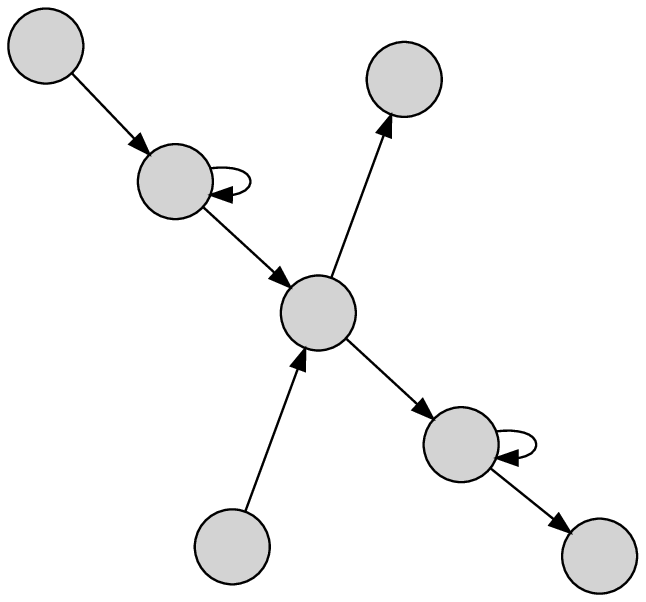}
\includegraphics[scale=0.25]{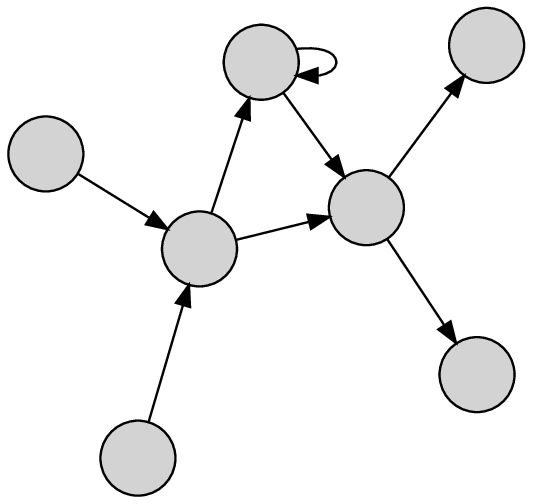}
\includegraphics[scale=0.25]{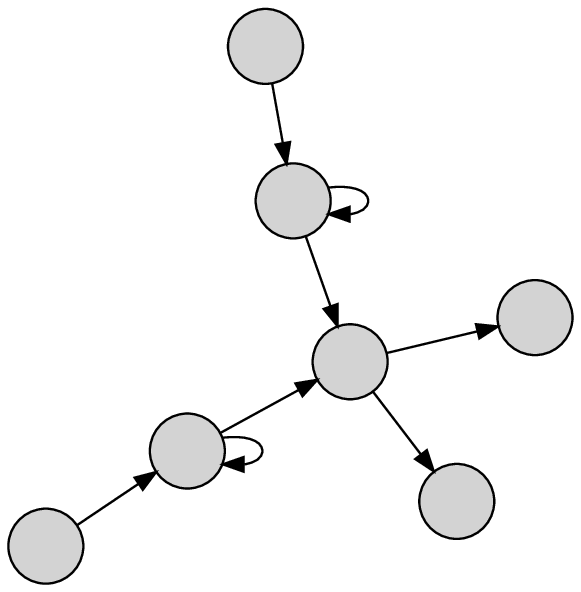}
\includegraphics[scale=0.25]{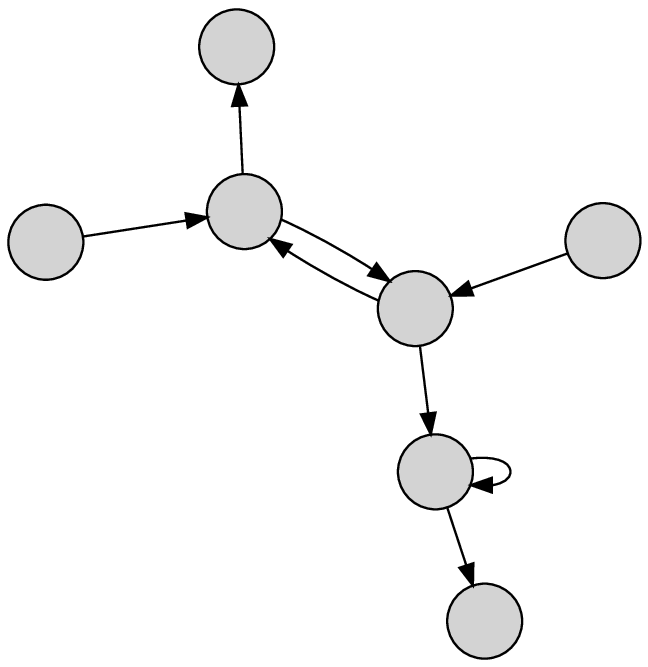}
\includegraphics[scale=0.25]{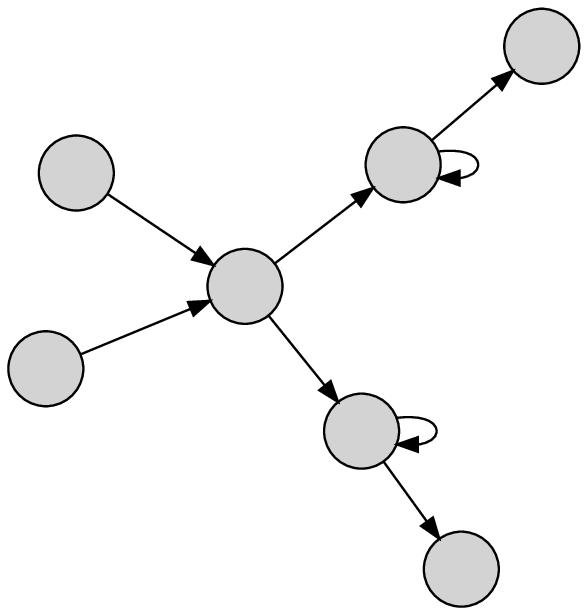}
\includegraphics[scale=0.25]{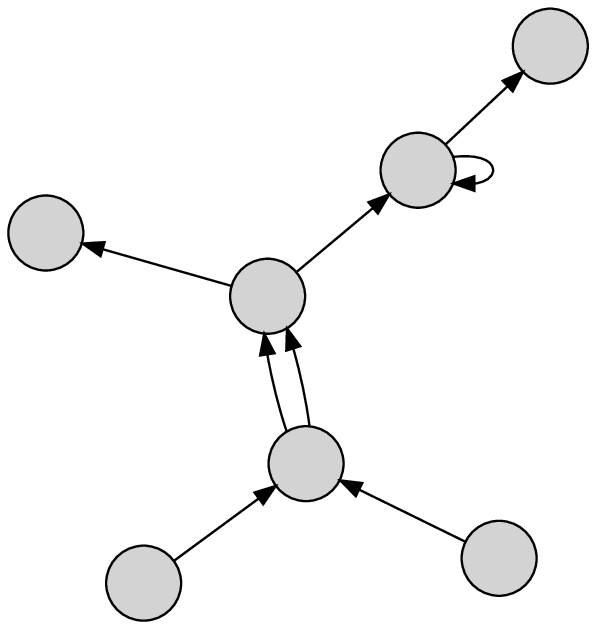}
\includegraphics[scale=0.25]{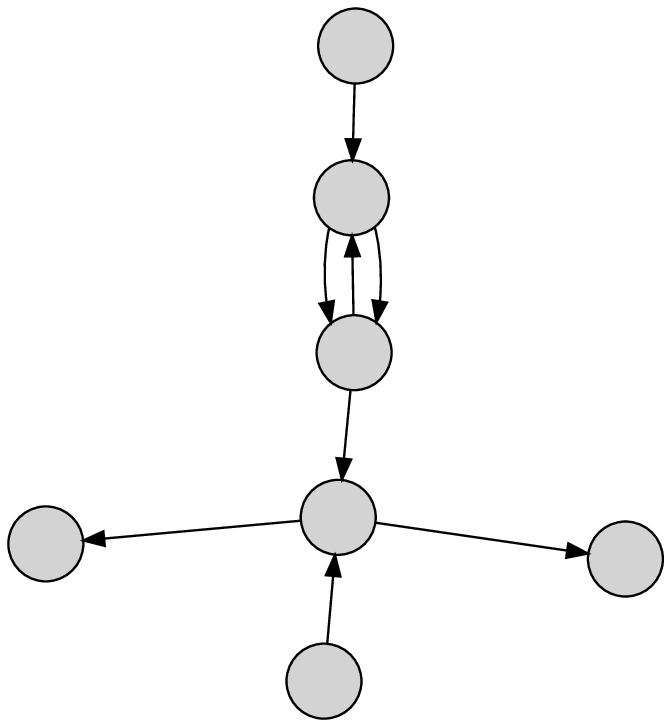}
\includegraphics[scale=0.25]{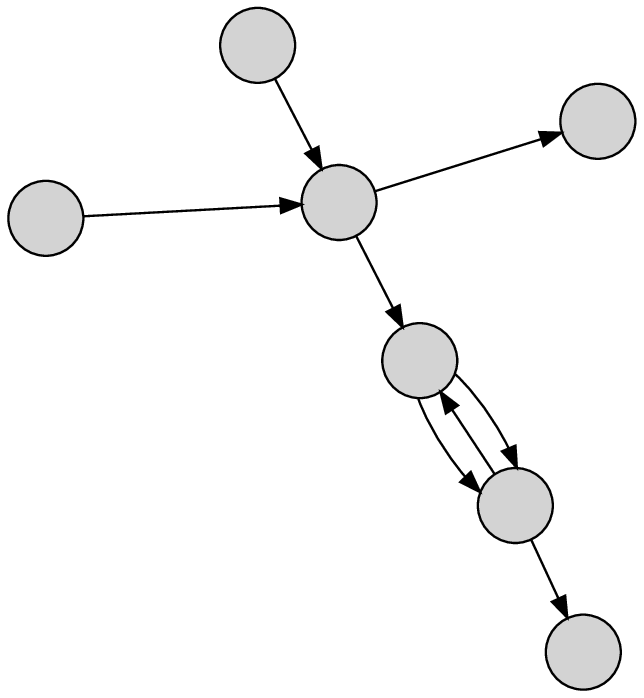}
\includegraphics[scale=0.25]{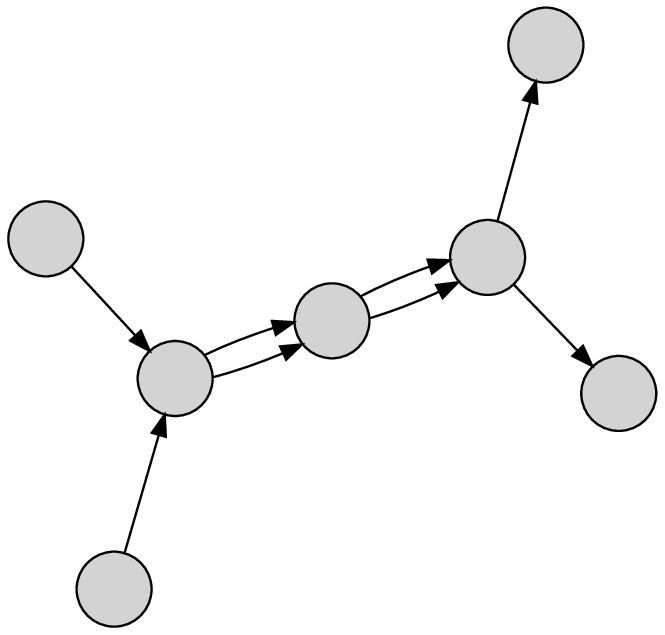}
\includegraphics[scale=0.25]{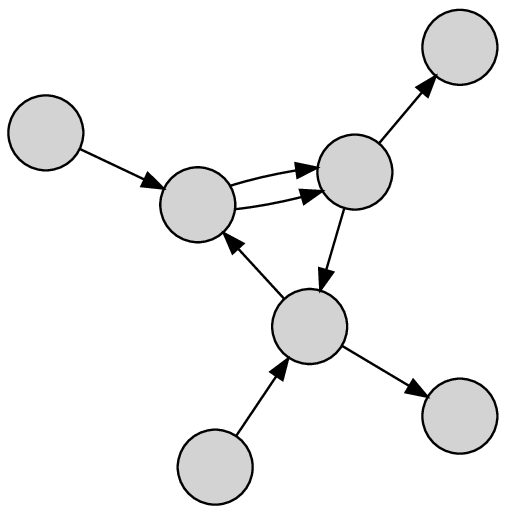}
\includegraphics[scale=0.25]{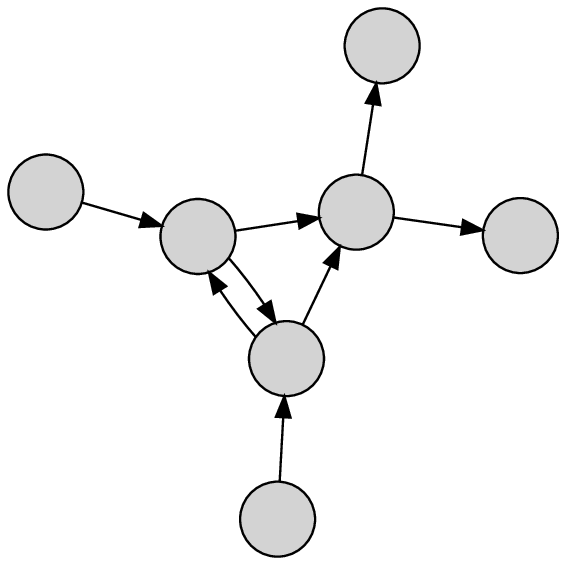}
\includegraphics[scale=0.25]{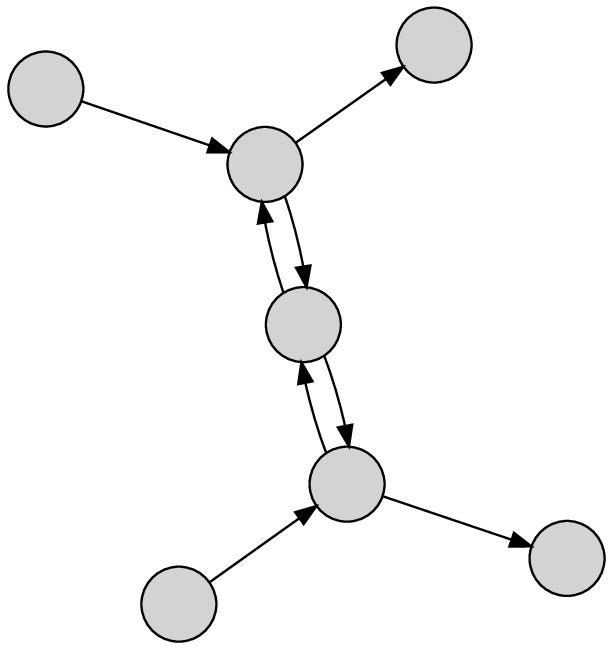}
\includegraphics[scale=0.25]{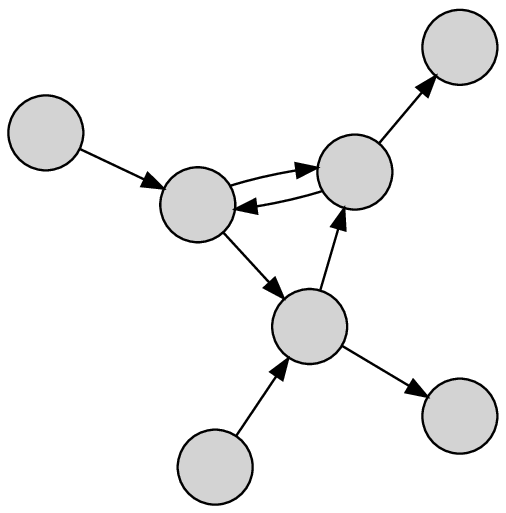}
\includegraphics[scale=0.25]{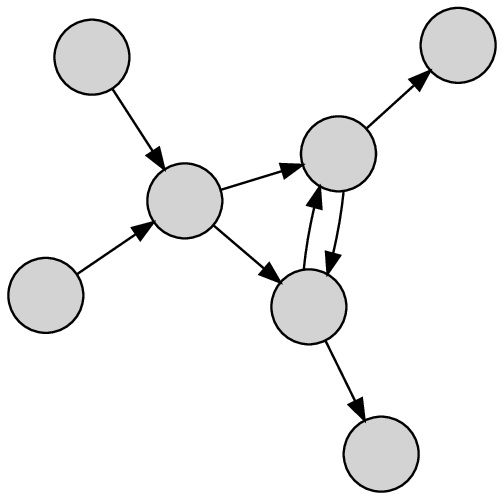}

A summary of the enumeration with two
arcs cut is in Table \ref{tab.F}.

\begin{table}[htb]
\begin{ruledtabular}
\begin{tabular}{r|rrr}
$n$ & $U_{2,+4}(n,1)$ & $L_{2,+4}(n,1)$ & $L_{2,+4}(n)$\\
\hline
4 & 0 & 0 & 0 \\
5 & 1 & 30 & 210 \\
6 & 4 & 1260 & 5220 \\
7 & 21  & 60480 & 180810 \\
\end{tabular}
\end{ruledtabular}
\caption{Connected unlabeled 2-regular digraphs after cutting two arcs, and 
corresponding labeled 2-regular digraphs (one or any number of components)
allowing loops and multiarcs.} 
\label{tab.F}
\end{table}

\section{Summary}
Table \ref{tab.main} summarizes the bare counts $U_2(n,c)$.

\begin{table}[htb]
\begin{ruledtabular}
\begin{tabular}{r|rrrrrrrrr|r}
$n\backslash c$ & $1$ & $2$& 3 & 4 & 5 & 6 & 7 & 8& 9& $U_2(n)$ \\
\hline
1 & 1 &&&&&&&&& 1\\
2 & 2 & 1 &&&&&&&& 3\\
3 & 5 & 2 & 1 &&&&&&& 8 \\
4 & 14 & 8 & 2 & 1 &&&&&& 25 \\
5 & 50 & 24 & 8 & 2 & 1 &&&&& 85\\
6 & 265 & 93 & 28 & 8 & 2 & 1 &&&& 397\\
7 & 1601 & 435 & 108 & 28 & 8 & 2 & 1 &&& 2183 \\
8 & 11984 & 2486 & 507 & 113 & 28 & 8 & 2 & 1&& 15129 \\
9 & 101884 & 17211 & 2811 & 527 & 113 & 28& 8 & 2 & 1& 122585 \\
\end{tabular}
\end{ruledtabular}
\caption{Unlabeled 2-regular digraphs $U_2(n,c)$ with $n$ nodes and $c$ (weak) components, allowing loops and multiarcs 
\cite[A306892,A006372]{sloane}.}
\label{tab.main}
\end{table}
The row sums 
1, 3, 8, 25, \ldots count the graphs with any number of components.
$U_2(4)$ and $U_2(5)$ have already been computed by Briggs
\cite{Briggsarxiv9703074,Briggsarxiv9607033}.
We observe that Briggs' extrapolations 
to
more than 5 nodes \cite{Briggsarxiv9808050}
underestimate the true number of graphs for 6 -- 9 nodes.

\appendix
\section{Machine Readable Tables}

The ancillary directory contains the information of the unlabeled
2-regular digraphs in files named \texttt{Reg}$n$\texttt{.txt.zip}, where
$n$ is the number of nodes. 
After unzipping, each file contains up to three successive
lines per graph:
\begin{enumerate}
\item
An arc list for $2n$ arcs, referring to a representative of the labeled graphs created
by the $\cal A$-group, with labels from $0$ up to $n-1$, in square brackets.
Each bracket contains a pair of numbers, separated by a comma; the first is
the tail node and second the head node of the arc. Multiarcs are rendered by repeating
such pairs.
\item
A capital \texttt{V}, the number of multiarcs in the graph (i.e., the number of entries larger
than 1 in the Adjacency Matrix), a blank,
the trace of the Adjacency Matrix (i.e., the number of loops), a blank, and the P\'olya Cycle
Index (a multinomial in the free variables $t_i$).
\end{enumerate}

Filtering the lines that start with \texttt{V0} in these files we obtain
1, 3, 8, 27,\ldots unlabeled, not necessarily connected, 2-regular digraphs with $n\ge 2$ nodes without multiarcs \cite[A005641]{sloane}.

Filtering the lines that start with \texttt{V0\textvisiblespace 0} in these files we obtain
1, 2, 5, 23,\ldots unlabeled, loopless, not necessarily connected, 2-regular digraphs with $n\ge 3$ nodes without multiedges \cite[A219889]{sloane}.

As a further check, filtering the lines that contain \texttt{\textvisiblespace 0} we obtain 1, 2, 6, 15, 68,\ldots graphs
on $n\ge 2$ nodes without loops (which may have multiarcs) \cite[A307180]{sloane}.

One application of this information yields the $r$-rooted unlabeled
2-regular digraphs by defining the generating function $r(x)=1+x$ for
the number of ways of labeling a node as 0 (not marked) or 1 (marked),
and then substituting $t_i\to r(x^i)$ in the cycle indices. Summation over
the cycle indices of all graphs of fixed $n$ generates Table \ref{tab.root}.

\begin{table}[htb]
\begin{ruledtabular}
\begin{tabular}{r|rrrrrrrrr}
$n\backslash r$ & $0$ & $1$ & $2$& 3 & 4 & 5 & 6 \\
\hline
1 & 1 & 1\\
2 & 3 & 3 & 3 \\
3 & 8 & 13 & 13 & 8\\
4 & 25 & 58 & 88 & 58 & 25\\ 
5 & 85 & 310 & 588 & 588 & 310 & 85 \\
6 & 397 & 1909 & 4626 & 6035 & 4626 & 1909 & 397 \\
7 & 2183 & 13843 & 40417 & 66471 & 66471 & 40417 & 13843 \\
8 & 15129 & 114821 & 395324 & 782257 & 975715 & 782257 & 395324 \\
\end{tabular}
\end{ruledtabular}
\caption{Unlabeled 2-regular digraphs $U_2^{(r)}(n)=U_2^{(n-r)}(n)$ 
with $n$ nodes and $0\le r\le n$ rooted 
nodes, allowing loops and multiarcs. $U_2^{(0)}(n) = U_2(n)$.}
\label{tab.root}
\end{table}

Another information drawn from the cycle indices are the parities
of the permutations of the node-labeled 2-regular digraphs that
leave the unlabeled graph intact. A permutation is odd, iff
there is an odd number of even-length cycles, i.e., iff for even index 
$i$ there is an odd exponent $j$ in the factor $t_i^j$ of the cycle index.
For products of different $i$ the product rule for the signs of cycles rules.
Example: one of the 25 graphs on 4 nodes has the cycle index $(t_1^4 +2t_1^2t_2 +3t_2^2 +2t_4)/8$.
\begin{itemize}
\item
The term $t_1^4$ is one permutation with no even index, so this is an even permutation.
\item
The term $2t_1^2t_2$ is two permutations with an even index with exponent 1, so this
represents two odd permutations.
\item
The term $3t_2^2$ is three permutations with an even index with exponent 2, so this
represents three even permutations.
\item
The term $2t_4$ is two permutations with an even index with exponent 1, so this
represents two odd permutations.
\end{itemize}

The ancillary directory contains the information of the connected unlabeled
fairly 2-regular digraphs of Section \ref{sec.scut} in files named \texttt{Regf}$n$\texttt{.txt.zip}, where
$n$ is the number of nodes, including the 2 nodes with degree 1. 
The arc lists contain $2n-3$ arcs and the lines starting with $V$
are the cycle indices.

The ancillary directory contains the information of the connected unlabeled
fairly 2-regular digraphs of Section \ref{sec.2cut} in files named \texttt{RegF}$n$\texttt{.txt.zip}, where
$n$ is the number of nodes, including the 4 nodes with degree 1. 
The arc lists contain $2n-6$ arcs and the lines starting with $V$
are the cycle indices.

\section{Java Code}
The ancillary directory contains the source code of the Java classes
in the \texttt{org/nevec/rjm subdirectory}.
The \texttt{javac} command
in the \texttt{Makefile} compiles the code. The other lines
in the \texttt{Makefile} demonstrate how the counts for the unlabeled
digraphs are obtained by calling two \texttt{main} classes in succession.
The first is \texttt{Regul} which generates incidence matrices
for digraphs of mixed connectivity and of specified degree squences
and prints the corresponding
edge/arc lists which is redirected to an intermediate file.
\texttt{Unlbl2RegGraphHashSet}
then reads such an
edge/arc list and reduces these sets to unlabeled graphs edge
lists by running through all row-column permutations and leaving
only one representative in its standard output.
The two \texttt{man}-pages \texttt{Regul.1} and \texttt{Unlbl2RegGraphHashSet.1}
provide details of options and arguments of the two classes.

The examples in the \texttt{Makefile} are not optimized for speed.
They construct the $L_{2,+2}(n,c)$ for example with 2 nodes of
degree 1 and $n-2$ nodes of in- and outdegree 2. Just
two geometries occur: the 2 nodes are connected to the same node or
to different nodes. So one could gather all cases by running
the program twice forgetting the 2 leaf nodes, constructing
graphs with $n-2$ nodes:
(i) a run with 1 node of indegree 1 and outdegree 1 and n-3 nodes
indegree and outdegree 2, (ii) a run with 1 node of indegree 1 and outdegree 2,
1 node of outdegree 1 and indegree 2, and n-4 nodes of indegree and outdegree 2.
After the reduction to unlabeled graphs the 2 leaf nodes 
could be re-attached at the nodes which are uniquely characterized
by their underrun of degrees.

\bibliography{all}

\end{document}